
\def\LaTeX{%
  \let\Begin\begin
  \let\End\end
  \let\salta\relax
  \let\finqui\relax
  \let\futuro\relax}

\def\UK{\def\our{our}\let\sz s}
\def\USA{\def\our{or}\let\sz z}

\LaTeX

\USA


\salta

\documentclass[twoside,12pt]{article}
\setlength{\textheight}{24cm}
\setlength{\textwidth}{16cm}
\setlength{\oddsidemargin}{2mm}
\setlength{\evensidemargin}{2mm}
\setlength{\topmargin}{-15mm}
\parskip2mm


\usepackage{amsmath}
\usepackage{amsthm}
\usepackage{amssymb}
\usepackage[mathcal]{euscript}
\usepackage{cite}

\usepackage[usenames,dvipsnames]{color}
%
%
%
\def\gianni #1{{\color{blue}#1}}
\def\juerg #1{{\color{PineGreen}#1}}          
\def\pier #1{{\color{red}#1}}     
\def\revis #1{{\color{red}#1}}
%
%
\let\gianni\relax
\let\pier\relax
\let\juerg\relax
\let\revis\relax



\bibliographystyle{plain}


%

\finqui

\newcommand{\beq}{\begin{equation}}
\newcommand{\eeq}{\end{equation}}
\newcommand{\beqa}{\begin{equation}}
\newcommand{\eeqa}{\end{equation}}
\let\non\nonumber


\def\step #1 \par{\medskip\noindent{\bf #1.}\quad}


\def\aand{\quad\hbox{and}\quad}

\def\rhs{right-hand side}


\def\multibold #1{\def\arg{#1}%
  \ifx\arg\pto \let\next\relax
  \else
  \def\next{\expandafter
    \def\csname #1#1#1\endcsname{{\bf #1}}%
    \multibold}%
  \fi \next}

\def\pto{.}

\def\multical #1{\def\arg{#1}%
  \ifx\arg\pto \let\next\relax
  \else
  \def\next{\expandafter
    \def\csname cal#1\endcsname{{\cal #1}}%
    \multical}%
  \fi \next}


\def\multimathop #1 {\def\arg{#1}%
  \ifx\arg\pto \let\next\relax
  \else
  \def\next{\expandafter
    \def\csname #1\endcsname{\mathop{\rm #1}\nolimits}%
    \multimathop}%
  \fi \next}

\multibold
qwertyuiopasdfghjklzxcvbnmQWERTYUIOPASDFGHJKLZXCVBNM.

\multical
QWERTYUIOPASDFGHJKLZXCVBNM.

\multimathop
dist div dom meas sign supp .


\def\Accorpa #1#2 #3 {\gdef #1{\eqref{#2}--\eqref{#3}}%
  \wlog{}\wlog{\string #1 -> #2 - #3}\wlog{}}


\def\qed{\hfill $\square$}
\def\tint {\int_0^t}

\def\xinto{\int_\Omega}
\def\txinto{\int_0^t\!\!\!\int_\Omega}
\def\texinto{\int_0^T\!\!\!\int_\Omega}
\def\oma{\Omega}
\def\ginto{\int_\Gamma}
\def\tginto{\int_0^t\!\!\!\int_\Gamma}
\def\teginto{\int_0^T\!\!\!\int_\Gamma}
\def\lzht{L^2(0,t;H)}
\def\lzvt{L^2(0,t;V)}

\def\zetg{\zeta_\Gamma}

\def\pt{\partial_t}
\def\pn{\partial_{\bf n}}
\def\ds{\,{\rm d}s}
\def\dt{\,{\rm d}t}
\def\dx{\,{\rm d}x}

\def\ug{\,u_\Gamma}

\def\muh{\mu^h}
\def\rh{\rho^h}
\def\rgh{\rho^h_\Gamma}
\def\zeh{\zeta^h}
\def\zeg{\zeta_\Gamma^h}
\def\etah{\eta^h}
\def\bug{\bar u_\Gamma}
\def\br{\bar\rho}
\def\brg{\bar\rho_\Gamma}
\def\bm{\bar\mu}
\def\ygh{y^h_\Gamma}
\def\vp{\varphi}
\def\etan{\eta_n}
\def\zetan{\zeta_n}
\def\zetang{\zeta_{n_\Gamma}}

\def\checkmmode #1{\relax\ifmmode\hbox{#1}\else{#1}\fi}
\def\aeO{\checkmmode{a.\,e.\ in~$\,\Omega$}}
\def\aeQ{\checkmmode{a.\,e.\ in~$\,Q$}}
\def\aeS{\checkmmode{a.\,e.\ on~$\,\Sigma$}}
\def\aeG{\checkmmode{a.\,e.\ on~$\,\Gamma$}}

\def\limrho0{\lim_{\rho\searrow -1}}
\def\limrho1{\lim_{\rho\nearrow +1}}

\def\rz{{\mathbb{R}}}
\def\nz{{\mathbb{N}}}


\def\Lx #1{L^{#1}(\Omega)}
\def\Hx #1{H^{#1}(\Omega)}

\def\Ldue{\Lx 2}

\def\Huno{\Hx 1}
\def\Hdue{\Hx 2}

\def\Hg{H_\Gamma}
\def\Vg{V_\Gamma}


\def\CQ {C^0(\overline Q)}
\def\CS {C^0(\overline\Sigma)}
\def\sol{{\cal S}}
\def\CJ{{\cal J}}
\def\TT{{\cal T}_{\tau}}
\def\TTN{{\cal T}_{\tau_N}}
\def\taun{\tau_N}


\let\TeXchi\chi                         
\newbox\chibox
\setbox0 \hbox{\mathsurround0pt $\TeXchi$}
\setbox\chibox \hbox{\raise\dp0 \box 0 }
\def\chi{\copy\chibox}


\def\r0g{\rho_{0_{|\Gamma}}}
\def\rg{\rho_\Gamma}
\def\pig{\pi_\Gamma}
\def\pigs{\pi_\Gamma'}

\def\ng{\nabla_\Gamma}
\def\delg{\Delta_\Gamma}
\def\dg{\,{\rm d}\Gamma}
\def\fg{f_\Gamma}
\def\fgas{f_\Gamma'}
\def\fgasz{f_\Gamma''}
\def\rtil{\tilde\rho}
\def\mtil{\tilde\mu}
\def\rgtil{\tilde\rho_\Gamma}
\def\mqtil{\tilde\mu_Q}
\def\rqtil{\tilde\rho_Q}
\def\rstil{\tilde\rho_\Sigma}
\def\zg{z_\Gamma}
\def\onSigma{_{|\Sigma}}

\def\coeff{1+2g(\rho)}

\def\uad{{\cal U}_{\rm ad}}

\Begin{document}


\title{\bf Optimal boundary control of a nonstandard\\ 
viscous Cahn-Hilliard system\\ with dynamic boundary 
condition\footnote{\revis{This work received a partial support from 
the MIUR-PRIN Grant 2015PA5MP7 ``Calculus of Variations'', 
the GNAMPA (Gruppo Nazionale per l'Analisi
Matematica, la Probabilit\`{a} e loro Applicazioni) of INDAM (Istituto Nazionale
di Alta Matematica) and the \pier{IMATI -- C.N.R. Pavia} for PC and GG.}}}

\author{}
\date{}
\maketitle
\begin{center}
\vskip-2cm
{\large\bf Pierluigi Colli$^{(1)}$}\\
{\normalsize e-mail: {\tt pierluigi.colli@unipv.it}}\\[.25cm]
{\large\bf Gianni Gilardi$^{(1)}$}\\
{\normalsize e-mail: {\tt gianni.gilardi@unipv.it}}\\[.25cm]
{\large\bf J\"urgen Sprekels$^{(2)}$}\\
{\normalsize e-mail: {\tt sprekels@wias-berlin.de}}\\[.45cm]
$^{(1)}$
{\small Dipartimento di Matematica ``F. Casorati'', Universit\`a di Pavia}\\
{\small via Ferrata \pier{5}, 27100 Pavia, Italy}\\[.2cm]
$^{(2)}$
{\small Weierstrass Institute for Applied Analysis and Stochastics}\\
{\small Mohrenstra\ss e\ 39, 10117 Berlin, Germany, and}\\
{\small Department of Mathematics, Humboldt-Universit\"at zu Berlin}\\
{\small Unter den Linden 6, 10099 Berlin. Germany}\\[.8cm]
\end{center}


\Begin{abstract}\noindent
In this paper, we study an optimal boundary control problem for a model for phase 
separation taking place in a
spatial domain that was introduced 
\gianni{by Podio-Guidugli in Ric.\ Mat.\ {\bf 55} (2006), pp.~105--118}. 
The model consists of a strongly  coupled system of nonlinear
parabolic differential equations, in  which products between  
the unknown functions and their time derivatives occur that are difficult to handle
analytically. In contrast to the existing control literature about this PDE system, we 
consider here a dynamic boundary condition involving
the Laplace--Beltrami operator for the order parameter of the system, which 
 models an additional nonconserving phase
transition occurring on the surface of the domain. We show the Fr\'echet differentiability
of the associated control-to-state operator in appropriate Banach spaces and derive results
on the existence of optimal controls and on first-order necessary optimality conditions in terms of a variational inequality and the adjoint state system.
\\[10mm]
\End{abstract}


\salta

\pagestyle{myheadings}
\newcommand\testopari{\sc Colli \ --- \ Gilardi \ --- \ Sprekels}
\newcommand\testodispari{\sc Optimal dynamic boundary control of a 
Cahn--Hilliard system}
\markboth{\revis{\testopari}}{\revis{\testodispari}}

\finqui


\section{Introduction}
\label{Intro}
\setcounter{equation}{0}
Let $\oma\subset\rz^3$ be an open, bounded and connected set with a \revis{smooth boundary $\Gamma$
and let} $Q:=\oma\times (0,T)$ and $\Sigma:=\Gamma\times (0,T)$. We denote by $\pn$, $\ng$,
$\delg$, the outward normal derivative, the tangential gradient, and the Laplace--Beltrami
operator on $\Gamma$, in this order. \pier{We consider the following}
optimal boundary control problem:

\vspace{3mm}\noindent
{\bf (CP)} \quad Minimize the (tracking-type) cost functional
\begin{align}
\label{cost}
&\CJ((\mu,\rho,\rg),\ug)\non\\
&:=\,\frac{\beta_1}2\,\|\mu-\hat\mu_Q\|^2_{L^2(Q)}\,+\,\frac{\beta_2}2\,
\|\rho-\hat\rho_Q\|^2_{L^2(Q)}\,+\,\frac{\beta_3}2\,\|\rg-\hat\rho_\Sigma\|^2_{L^2(\Sigma)}
\quad\non\\
&+\frac{\beta_4}2\,\|\rho(T)-\hat\rho_\Omega\|^2_{L^2(\oma)}\,+\,\frac{\beta_5}2
\,\|\rg(T)-\hat\rho_\Gamma\|^2_{L^2(\Gamma)}\,+\,\frac{\beta_6}2\,\|\ug\|^2_{L^2(\Sigma)}\quad
\end{align}
over a suitable set  $\uad\subset (H^1(0,T;L^2(\Gamma))\cap 
L^\infty(\Sigma))$ of admissible controls $\ug$ (to be specified later),
subject to the state system
\begin{eqnarray}
  && \bigl(\coeff \bigr) \, \pt\mu
  + \mu \, g'(\rho) \, \pt\rho
  - \Delta\mu =0
  \gianni{\aand \mu\geq 0}
  \quad\mbox{in $Q$,}
  \label{ss1}
  \\
	&&\pn\mu =0 \quad\mbox{on $\Sigma$},\label{ss2}\\
  && \pt\rho - \Delta\rho + f'(\rho) +\pi(\rho) = \mu \, g'(\rho)	
	\quad\mbox{in $Q$},
  \label{ss3}
  \\
\label{ss4}  
  && \pn\rho+\pt\rg+\fgas(\rg)+
	\pig(\rg) -\delg	\rg= \ug,
	\quad \rg=\rho\onSigma,\quad \mbox{on $\Sigma$},\\
  \label{ss5}
    && \mu(0) = \mu_0,\quad \rho(0)=\rho_0,
	\quad\mbox{in $\Omega$}, \quad
    \rho_\Gamma(0) = \rho_{0_{|\Gamma}}\quad\mbox{on $\Gamma$.}
 \end{eqnarray}
Here, $\beta_i$, $1\le i\le 6$, are nonnegative weights, and $\,\hat\mu_Q, \hat\rho_Q\in L^2(Q)$,
$\,\hat\rho_\Sigma\in L^2(\Sigma)$, $\,\hat\rho_\Omega\in L^2(\oma)$, 
and $\,\hat\rho_\Gamma\in L^2(\Gamma)$\, are prescribed target functions. 
Although more general cost functionals could be admitted for large parts of the subsequent
analysis, we restrict ourselves to the above situation for the sake of a simpler exposition.

The physical background behind the control problem {\bf (CP)} is the following: the state system \eqref{ss1}--\eqref{ss5} constitutes a model for phase separation taking place in the container $\oma$ and originally introduced in \cite{PG}. In this
connection, the unknowns $\,\mu\,$ and $\,\rho\,$ denote the associated chemical potential, which 
in this particular model has to be nonnegative \gianni{(see~\eqref{ss1})}, and the order 
parameter of the phase separation process, which is usually the volumetric density of one of the involved phases.
We assume that $\,\rho\,$ is normalized in such a way as to attain its values in the interval $(-1,1)$. 
The nonlinearities $\pi,\pig,g$ are assumed to be smooth in $[-1,1]$, while $\,f\,$ and $\,f_\Gamma\,$ 
are \revis{convex} potentials defined in $(-1,1)$, whose derivatives 
$f',\fgas$ are singular at the endpoints $r=-1$ and $r=1$. A typical case is given by the \emph{logarithmic
potential}
\beq\label{logpot}
f(r)=f_\Gamma(r)=\hat c\,((1+r)\,\log(1+\gianni r)+(1-r)\,\log(1-r)), \,\mbox{ with 
a constant $\,\hat c>0$}.
\eeq
The state system (\ref{ss1})--(\ref{ss5}) is singular, with highly nonlinear
and nonstandard coupling. In particular, unpleasant nonlinear terms involving 
time derivatives occur in (\ref{ss1}),
and the expression\gianni s $f'(\rho)$ and $\fgas(\rg)$  in (\ref{ss3}), \eqref{ss4}  may become singular. \revis{We point out that in the applications the sum of $f$ and of an antiderivative of $\pi$ usually represents a double-well or multi-well potential, and the same can be said for the analogous functions on the boundary.}

The state system has been the subject of intensive study in the past years for the case that \eqref{ss4} is 
replaced by a zero Neumann condition. In this \pier{connection}, we refer the reader to  
\cite{CGKPS,CGKS1,CGKS2,CGPS3, CGPS6, CGPS7, CGPS4,CGPS5}. In \cite{CGPSco} an associated control
problem with a distributed control  in  \eqref{ss1} was investigated for the special case $g(\rho)=\rho$, 
and in~\cite{CGSco1} the corresponding case of a boundary control in \eqref{ss2} was studied. A nonlocal
version, in which the Laplacian $\,-\Delta\rho\,$ in \eqref{ss3} was replaced by a nonlocal operator, was
discussed in the recent contributions \cite{CGS3, CGS4,CGS4neu}. 

In all of the works cited above a zero Neumann condition was assumed for the order parameter $\,\rho$.
In contrast to this, we study in this paper the case of the dynamic boundary condition \eqref{ss4}. 
It models a nonconserving phase transition taking place on the boundary, which could be, e.\,g., induced
by an interaction between bulk and wall. The associated total free energy of the phase separation process is  
the sum of a bulk and a surface energy and has the form
\begin{align}\label{ftot}
&{\cal F}_{\rm tot}[\pier{{}\mu(t),{}}\rho(t),\rg(t)] \non\\[1mm] 
&:=\xinto \Big(f(\rho(x,t))+\hat\pi(\rho(x,t))\,\pier{{}-\,\mu(x,t)\,g(\rho(x,t))}\,+\,\frac 12|\nabla\rho(x,t)|^2
\Big)\dx\non\\[1mm] 
&+\ginto\Big(f_\Gamma(\rg(x,t))+\hat\pi_\Gamma(\rg(x,t))\,\pier{{}-u_\Gamma(x,t)\,\rg(x,t)}\,+\,\frac 12|\nabla_\Gamma\rg(x,t)|^2
\Big)\dg\,,
\end{align}
for $t\in [0,T]$, where $\hat\pi(r)=\int_0^r\pi(\xi){\rm d}\xi\,$ and \,$\hat\pi_\Gamma
(r)=\int_0^r\pig(\xi){\rm d}\xi$. 

In the recent contribution \cite{CGSneu}, the state system \eqref{ss1}--\eqref{ss5} was studied
systematically concerning existence, uniqueness, and regularity. Notice that in \cite{CGSneu}
more general nonlinearities were admitted, including the case that $f,f_\Gamma$ could be
nondifferentiable indicator functions (in which case $f'(\rho)$ and $\fgas(\rg)$ have to be interpreted as elements
of the (possibly multivalued) subdifferentials of $\,f\,$ at $\,\rho\,$ and of $\,f_\Gamma\,$
at $\,\rg$, respectively, so that \eqref{ss3} and \eqref{ss4} have to be understood as 
differential inclusions).    

The mathematical literature on control problems for phase field systems involving equations
of viscous or nonviscous Cahn--Hilliard type is still scarce and quite recent. We refer in this connection to the works \cite{CFGS1, CFGS2,CGS1,CGS2,HW,wn99}. Control problems
for convective Cahn--Hilliard systems were studied in\pier{\cite{RS,ZL1,ZL2}}, and a few
analytical contributions were made to the coupled Cahn--Hilliard/Navier--Stokes system
(cf. \cite{FRS,HW3, HW1,HW2}). The contribution \cite{CGRS}
dealt with the optimal control of a  Cahn--Hilliard type system arising in the modeling of
solid tumor growth. For the optimal control of Allen--Cahn equations with dynamic
boundary condition\pier{s}, we refer to \cite{CFS,CS}.

\vspace{2mm}
The paper is organized as follows: in Section 2, we formulate the relevant assumptions on the
data of the control problem {\bf (CP)}, and we prove a strong stability result for the
state system \eqref{ss1}--\eqref{ss5}. In Section 3, we prove the Fr\'echet differentiability of  
the control-to-state operator in appropriate Banach spaces. Section 4 then brings the main
results of this paper, namely, the existence of optimal controls and the derivation of
the first-order necessary conditions of optimality. \revis{In our approach to the optimal control problem, we are of course inspired by the general theory contained in the monography~\cite{Tr}.}

Throughout the paper, we denote for a general Banach space $\,X\,$ by $\|\cdot\|_X$ its norm and by $X'$ its dual space.
The only exemption from this convention are the norms of the $L^p$ spaces and of  their powers, which we often denote by
$\|\cdot\|_p$, for $1\le p\le +\infty$.   
Moreover, we repeatedly utilize the continuity of the embedding $\Huno \subset\Lx p$
for $1\leq p\leq 6$ and the related Sobolev inequality
\beq
  \|v\|_p \leq \gianni{C_\Omega} \|v\|_{\Huno}
  \quad \hbox{for every $v\in \Huno$ and $1\leq p \leq 6$,}
  \label{sobolev}
\eeq
where $\gianni{C_\Omega}$ depends only on~$\Omega$. Notice that these embeddings  
are compact for $1\le p<6$. We also recall that
the embedding \,$\Hdue\subset C^0(\overline{\oma})$ is compact.
Furthermore, we make repeated use of \pier{H\"older's inequality and of the elementary Young inequality
\begin{align}
\label{Young}
& |ab| \leq \gamma\, |a|^2 + \mbox{$\frac 1{4\gamma}$}
 \, |b|^2\,\quad\mbox{for every $\,a,b\in \rz\,$ 
and $\,\gamma>0$},
	\end{align}
and we set
\begin{align}
  &Q_t := \Omega \times (0,t),\quad\Sigma_t:=\Gamma\times (0,t), 
  \quad \hbox{for $\,t\in (0,T].$}
  \label{defQt}
\end{align}
About time derivatives of a time-dependent function $v$, we warn the reader 
that we will use both the notations  $\pt v, \, \pt^2 v $ and the shorter 
ones $v_t, \, v_{tt} $.}

\section{General assumptions and results\\ for the state system}
\setcounter{equation}{0}

In this section, we formulate the general assumptions for the data of the control problem
{\bf (CP)}, and we state some preparatory results for the state system 
\eqref{ss1}--\eqref{ss5}. To begin with, we
introduce some denotations. We set
\begin{align*}
&H:=\Ldue,\quad V:=\Huno, \quad W:=\{w\in\Hdue:\mbox{\,$\pn w=0$\, on $\,\Gamma$}\},\\
&\Hg:=L^2(\Gamma),\quad\Vg:=H^1(\Gamma),\quad {\cal V}:=\{v\in V:v_{|\Gamma}\in \Vg\},
\end{align*}
and endow these spaces with their standard norms. Notice that we have $V\subset H\subset V'$ and 
$\Vg\subset\Hg\subset \Vg'$ with dense, continuous and compact embeddings.

We make the following general assumptions:

\vspace{3mm}\noindent
(A1) \quad $\mu_0\in W,\,\,\,\mu_0\ge 0$ \,a.\,e. in\, $\oma$,
\,\,$\,\rho_0\in \Hdue$, \,\,\,$\rho_{0_\Gamma}:=\rho_{0_{|\Gamma}}\in H^2(\Gamma)$, and
\beq\label{initial}
-1\,<\,\min_{x\in \overline\oma}\,\rho_0(x), \quad \max_{x \in\overline\oma}\rho_0(x)\,<\,+1.
\eeq

\vspace{1mm}\noindent
(A2) \quad $\pi,\pig\in C^2[-1,+1]$\pier{;} $\,g\in C^3[-1,+1]\,$ is nonnegative and concave
on $[-1,+1]$.

\vspace{1mm}\noindent
(A3) \quad  $f,\fg\in C^3(-1,+1)$\, are nonnegative and convex, satisfy \,$f(0)=\fg(0)=0$,
and \\
\hspace*{14mm}there are constants $\delta>0$ and $C_\Gamma\ge 0$ such that 
\beq \label{domina1} 
|f'(r)|\,\le\,\delta\,|\fgas(r)|\,+\,C_\Gamma \quad\forall\,r\in (-1,+1). 
\eeq
\hspace*{15mm}Moreover, it holds that
\begin{align}
\label{singul}
\lim_{r\searrow -1}\,f'(r)=\lim_{r\searrow -1}\,\fgas(r)=-\infty, \quad
\lim_{r\nearrow +1}\,f'(r)=\lim_{r\nearrow +1}\,\fgas(r)=+\infty.
\end{align}

\vspace{1mm}\noindent 
(A4) \quad $\uad=\Bigl\{\ug\in H^1(0,T;\Hg)\cap L^\infty(\Sigma):\,
\mbox{$u_*\le \ug\le u^*$ \,a.\,e. on 
$\Gamma$ \,and}$\\
\pier{\hspace*{9cm}$ \|{\ug}\|_{H^1(0,T;\Hg)\cap L^\infty(\Sigma)}\,\le R_0\Bigr\}$,}\\[2mm]
\hspace*{14mm}where $\,u_*,u^*\in L^\infty(\Sigma)\,$ and $R_0>0$  are such that 
$\uad\not=\emptyset$.

\vspace{1mm}\noindent
(A5) \quad Let $R>0$ be a constant such that
$\uad\subset {\cal U}_R$ with  the open ball\\[1mm]
\hspace*{12mm} $\,{\cal U}_R
:=\left\{\revis{\ug} \in H^1(0,T;\Hg)\cap L^\infty(\Sigma):\|\revis{\ug}\|_{H^1(0,T;\Hg)\cap L^\infty(\Sigma)}
\,<\,R\right\}$.  

\vspace{1mm}\noindent
(A6) \quad The constants $\beta_i$, $1\le i\le \revis{6}$, are nonnegative, and we have that
$\,\hat\mu_Q,\hat\rho_Q\in L^2(Q)$, \\
\hspace*{14.5mm}$\hat\rho_\Sigma\in L^2(\Sigma)$, $\hat\rho_\oma
\in L^2(\oma)$, and $\,\hat\rho_\Gamma\in L^2(\Gamma)$.

\vspace{5mm} The assumption (A5) is rather a denotation. We also remark that (A3) entails, in particular, 
that $\,f'(0)=f_\Gamma'(0)=0$, and it is easily seen that (A3) is fulfilled for the 
logarithmic potential\gianni s~\eqref{logpot}, 
\gianni{even with different $\hat c$'s for $f$ and~$f_\Gamma$}. 
In addition, if we assume that $\ug\in{\cal U}_R$, 
then it follows from the assumptions (A1)--(A3)
that \cite[Thm.~2.4]{CGSneu} can be applied. In fact, a closer inspection of the
proof of \cite[Thm.~2.4]{CGSneu} reveals that the following result holds true.

\vspace{5mm}\noindent
{\sc Theorem 2.1:} \quad\,\,{\em Suppose that} (A1)--(A5) {\em are fulfilled. Then the state
system} \eqref{ss1}--\eqref{ss5} {\em has for every $\ug\in {\cal U}_R$ a unique solution triple 
$\,(\mu,\rho,\rg)\,$ such that}
\begin{align}
\label{regmu}
&\mu\in H^1(0,T;H)\cap C^0([0,T];V)\cap L^2(0,T;W)\cap L^\infty(Q),\\
\label{regrho}
&\rho\in W^{1,\infty}(0,T;H)\cap H^1(0,T;V)\cap L^\infty(0,T;H^2(\oma)),\\
\label{regrg}
&\rg\in W^{1,\infty}(0,T;\Hg)\cap H^1(0,T;\Vg)\cap L^\infty(0,T;H^2(\Gamma)).
\end{align}
{\em Moreover, there is a constant $K_1^*>0$, 
which depends only on $R$ and the data of the system, such that}
\begin{align}\label{ssbounds1}
&\|\mu\|_{H^1(0,T;H)\cap C^0([0,T];V)\cap L^2(0,T;W)\cap L^\infty(Q)}\,
+\,\|\rho\|_{ W^{1,\infty}(0,T;\pier{H})\cap H^1(0,T;V)
\cap L^\infty(0,T;H^2(\oma))}\non\\
&+\,\|\rg\|_{W^{1,\infty}(0,T;\Hg)\cap H^1(0,T;\Vg)\cap L^\infty(0,T;H^2(\Gamma))}
\,\le\,K^*_1\,,
\end{align}
{\em for every solution triple $(\mu,\rho,\rg)$ corresponding to some $\ug\in{\cal U}_R$. In addition,
there are constants $r_*,r^*$, which depend only on $R$ and the data of the system,
such that}
\beq\label{separ}
-1<r_*\,\le\,\rho(x,t)\le r^*<+1 \quad\mbox{for every \,$(x,t)\in \overline Q$},
\eeq
{\em for every solution triple $(\mu,\rho,\rg)$ corresponding to some $\ug\in{\cal U}_R$.}

\vspace{3mm}
\noindent
{\sc Remark 2.2:} \quad\,\,It follows from well-known embedding results (cf.
\cite[Sect.~8,~Cor.~4]{Simon}) that $(H^1(0,T;V)\cap L^\infty(0,T;H^2(\oma)))
\subset C^0([0,T];H^s(\oma))$, for $0<s<2$. Therefore, $\rho\in\CQ$, and
thus $\rg\in \CS$. Moreover, there is a constant $K_2^*>0$, which again depends
only on $R$ and the data, such that
\begin{align}\label{ssbounds2}
&\max_{0\le i\le 3}\left(\bigl\|f^{(i)}(\rho)\bigr\|_{\CQ}\,+\,
\bigl\|f_\Gamma^{(i)}(\rg)\bigr\|_{\CS}\right)\,+\,\max_{0\le i\le 3}\,\bigl\|
g^{(i)}(\rho)\bigr\|_{\CQ}\non\\
&+\,\max_{0\le i\le 2}\left(\bigl\|\pi^{(i)}(\rho)\bigr\|_{\CQ}\,+\,\bigl\|\pi^{(i)}
_\Gamma(\rg)\bigr\|_{\CS}\right)\le\,K^*_2\,,
\end{align}
for every solution triple $(\mu,\rho,\rg)$ corresponding to some $\ug\in{\cal U}_R$.
In addition, we have \pier{that} $\,1+2g(\rho)\in\CQ$, where $\,1\le 1+2g(\rho)\le
1+2\|g(\rho)\|_{\CQ}\,$ on $\,\overline Q$. Hence, rewriting \eqref{ss1} \pier{as}
$$
\pt\mu - \mbox{$\frac 1{1+2g(\rho)}$}\,\Delta\mu=z \quad\pier{\hbox{in } Q,}
$$
where it is easily seen that 
$\,z:=-(1+2g(\rho))^{-1}\,\mu\,g'(\rho)\,\pt\rho\in L^\infty(0,T;H)\cap
L^2(0,T;L^6(\oma))$, we may thus infer from \cite[\revis{Thm.~2.1}]{DHP} the additional regularity
\begin{align}
\label{addregmu}
\mu\in W^{1,p}(0,T;H)\cap H^1(0,T;L^6(\oma))\cap L^p(0,T;W)\cap L^2(0,T;W^{2,6}
(\oma))\quad\non\\ \pier{\mbox{for every }\,p\in [1,+\infty).}
\end{align}
Moreover, denoting by 
\beq
  \gianni{{\cal X}:=H^1(0,T;\Hg)\cap L^\infty(\Sigma)}
  \label{defX}
\eeq
the control space for the remainder of this paper, we conclude from Theorem~2.1 that
the control-to-state operator $\,{\cal S}:\ug\mapsto (\mu, \rho,\rg)$, 
where it is understood that $\rg=\rho\onSigma$ \pier{on~$\Sigma$}, is a 
well-defined mapping between $\,{\cal U}_R\subset {\cal X}\,$ and the space specified by
the regularity properties \gianni{\eqref{regmu}--\eqref{regrg}}. 

\vspace{3mm}
For later use, we cite a known auxiliary result (cf. \cite[Thm.~2.2]{CS}).

\vspace{4mm}\noindent
{\sc Lemma 2.3:} \quad\,\,{\em Suppose that functions $y_0\in{\cal V}$, $a\in L^\infty(Q)$,
$a_\Gamma\in L^\infty(\Sigma)$, $\sigma\in L^2(Q)$ and $\sigma_\Gamma\in L^2(\Sigma)$ are 
given. Then the linear initial-boundary value problem}
\begin{align}
\label{CoSp1}
&\pt y-\Delta y\,+\,a\,y\,=\,\sigma \quad\mbox{{\em a.\,e. in }}\,Q,\\
\label{CoSp2}
&\pn y+\pt y_\Gamma-\delg y_\Gamma\,+\,a_\Gamma\,y_\Gamma\,=\,\sigma_\Gamma,\quad y_\Gamma
=y\onSigma, \quad\mbox{{\em a.\,e. on }}\,\Sigma,\\
\label{CoSp3}
&y(0)=y_0 \quad\mbox{{\em a.\,e. in }}\,\oma,\quad y_\Gamma(0)=y_{0|\Gamma}
\quad\mbox{{\em a.\,e. on }}\,\Gamma, 
\end{align}
{\em has a unique solution pair satisfying $\,y\in H^1(0,T;H)\cap C^0([0,T];V)
\cap L^2(0,T;H^2(\oma))\,$ and
$y_\Gamma\in H^1(0,T;\Hg)\cap C^0([0,T];\Vg)\cap L^2(0,T;H^2(\Gamma))$. Moreover, there is a constant $C_L>0$ such that 
the following holds true: whenever $y_0=0$ and $(y,y_\Gamma)$ is the corresponding
solution to} \eqref{CoSp1}--\eqref{CoSp3}, {\em then}
\begin{align}
\label{CoSp4}
&\|y\|_{H^1(0,t;H)\cap C^0([0,t];V)\cap L^2(0,t;H^2(\oma))}
\,+\,\|y_\Gamma\|_{H^1(0,t;\Hg)\cap C^0([0,t];\Vg)\cap L^2(0,t;H^2(\Gamma))}\non\\
&\le\,C_L\left(\|\sigma\|_{L^2(Q_t)}\,+\,\|\sigma_\Gamma\|_{L^2(\Sigma_t)}\right)\quad\forall\,t\in (0,T]\,.
\end{align}

\vspace{2mm} 
We are now going to investigate the stability properties of the state system 
\eqref{ss1}--\eqref{ss5}. We have the following
result.

\vspace{5mm}\noindent
{\sc Theorem 2.4:}  \quad\,\,{\em Suppose that the assumptions} (A1)--(A5) {\em are 
fulfilled, and assume that $u_{\Gamma_1},u_{\Gamma_2}\in {\cal U}_R$ are given and that
$\,(\mu_i,\rho_i,\rho_{i_\Gamma})={\cal S}(u_{\Gamma_i})$, $i=1,2$, are the corresponding 
unique solutions to} \eqref{ss1}--\eqref{ss5}. {\em Then we have for every
$t\in (0,T]$  the estimate}
\begin{align}
\label{stability}
&\|\mu_1-\mu_2\|_{H^1(0,t;H)\cap C^0([0,t];V)\cap L^2(0,T;W)}\,+\,
\|\rho_1-\rho_2\|_{H^1(0,t;H)\cap C^0([0,t];V)\cap L^2(0,t;\Hdue)}\non\\[1mm]
&+\,\|\rho_{1_\Gamma}-\rho_{2_\Gamma}\|_{H^1(0,t;\Hg)\cap C^0([0,t];\Vg)
\cap L^2(0,t;H^2(\Gamma))}
\,\le\,K^*_3\,\|u_{\Gamma_1}-u_{\Gamma_2}\|_{L^2(\pier{\Sigma_t})}\,,
\end{align}
{\em with a constant $K^*_3>0$ that depends only on $\,R\,$ and the data of the system.}

\vspace{3mm}
\noindent
{\sc Proof:} \quad\,Let $t \in (0,T]$ be fixed, and suppose that 
$u_{\Gamma_1},u_{\Gamma_2}\in {\cal U}_R$ are given and that
$\,(\mu_i,\rho_i,\rho_{i_\Gamma})={\cal S}(u_{\Gamma_i})$, $i=1,2$, are the corresponding 
unique solutions to \eqref{ss1}--\eqref{ss5} having the regularity properties 
\eqref{regmu}--\eqref{regrg} and \eqref{addregmu}. Observe that then the global
bounds \eqref{ssbounds1} and \eqref{ssbounds2}
hold true for both solutions. In the following, we will make repeated use of these bounds 
without further reference. We will also denote by $\,C>0\,$ constants that depend only on the data,
on $\,\|u_{\Gamma_i}\|_{\cal X}$, and on the norms of $(\mu_i,\rho_i,\rho_{i_\Gamma})$
in the spaces specified in \gianni{\eqref{regmu}--\eqref{regrg}} and \eqref{addregmu}. 
Now put 
\begin{align*}
&\mu:=\mu_1-\mu_2, \quad
\rho:=\rho_1-\rho_2, \quad \rho_\Gamma:=\rho_{1_\Gamma}-\rho_{2_\Gamma},\quad
\ug:=u_{\Gamma_1}-u_{\Gamma_2}, \\[1mm]
&\psi:=f'+\pi, \quad \psi_\Gamma:=f_\Gamma'+\pig.
\end{align*}
Then the following system is satisfied:
\begin{align}\label{diff1}
&(1+2g(\rho_1))\,\pt\mu\,+\,g'(\rho_1)\,\pt\rho_1\,\mu -\Delta\mu\,+\,2(g(\rho_1)-g(\rho_2))
\,\pt\mu_2\non\\
&\quad+\,\mu_2(g'(\rho_1)-g'(\rho_2))\,\pt\rho_1\,+\,\mu_2\,g'(\rho_2)\,\pt\rho\,=\,0
\quad\mbox{a.\,e. in } \,Q,\\[1mm]
\label{diff2}
&\pn\mu\,=\,0 \quad\mbox{a.\,e. on }\,\Sigma,\quad\mu(0)=0\quad\aeO,\\[1mm]
\label{diff3}
&\pt\rho-\Delta\rho=\psi(\rho_2)-\psi(\rho_1)+\mu\,g'(\rho_1)+\mu_2\,
(g'(\rho_1)-g'(\rho_2))\quad\aeQ,\\[1mm]
\label{diff4}
&\pn\rho+\pt\rg-\Delta_\Gamma\rg=\psi_\Gamma(\rho_{2_\Gamma})-\psi_\Gamma(\rho_{1_\Gamma})
+\ug, \quad \rg=\rho\onSigma, \quad\aeS,\\[1mm]
\label{diff5}
&\rho(0)=0 \quad\aeO,\quad \rg(0)=0 \quad\aeS.
\end{align}	
We will now prove a series of estimates in order to establish the validity of \eqref{stability}.
At first, we observe that
\begin{align}
\label{diffi1}
&\max_{0\le i\le 2}\left|g^{(i)}(\rho_1)-g^{(i)}(\rho_2)\right|\,+\,\max_{0\le i\le 1}\left
|\psi^{(i)}(\rho_1)-\psi^{(i)}(\rho_2)\right|\,\le \,C\,|\rho| \quad\aeQ,\\[1mm]
\label{diffi2}
&\max_{0\le i\le 1}\left
|\psi_\Gamma^{(i)}(\rho_{1_\Gamma})-\psi_\Gamma^{(i)}(\rho_{2_\Gamma})\right|
\,\le \,C\,|\rg| \quad\aeS.
\end{align}

\vspace{2mm}\noindent
The first estimate can be inferred from the stability result of \cite[Thm.~2.4]{CGSneu},
namely that
\begin{align}\label{stabu1}
&\|\mu\|_{L^\infty(0,t;H)\cap L^2(0,t;V)}\,+\,\|\rho\|_{H^1(0,T;H)\cap 
C^0([0,t];V)\cap L^2(0,t;\Hdue)}\non\\[1mm]
&+\,\|\rg\|_{H^1(0,t;\Hg)\cap C^0([0,t];\Vg)\cap L^2(0,t;H^2(\Gamma))}
\,\le\,C\,\|\ug\|_{L^2(0,t;\Hg)}\,.
\end{align}

Next, we \juerg{add $\mu$ to both sides of~\eqref{diff1}, then multiply} by $\pt\mu$ and
integrate over $Q_t$ to obtain that 
\begin{align}\label{stabu2}
&\txinto|\pt\mu|^2\dx\ds\,+\,\frac 12\,\juerg{\|\mu(t)\|_V^2}\,\le I_1+I_2+I_3+I_4,
\end{align}
where the quantities $I_j$, $1\le j\le 4$, will be specified and estimated below. At
first, we employ the continuity of the \juerg{embeddings $V\subset L^4(\oma)\subset L^2(\oma)$}, as well as
H\"older's and Young's inequalities, to conclude that
\begin{align}
\label{stabu3}
I_1:&=\txinto \juerg{\bigl( 1 - g'(\rho_1)\,\pt\rho_1 \bigr)}\,\mu\,\pt\mu\dx\ds\,
  \non\\
  &\le\,C\tint \juerg{\bigl( 1+\|\pt\rho_1(s)\|_4 \bigr)}
\,\|\mu(s)\|_4\,\|\pt\mu(s)\|_2\ds\non\\
&\le\,\frac 16\txinto|\pt\mu|^2\dx\ds\,+\,C\tint\juerg{\bigl( 1+\|\pt\rho_1(s)\|_V^2\bigr)}\,\|\mu(s)\|_V^2\ds\,. 
\end{align} 
Similarly, by also using \eqref{diffi1} and \eqref{stabu1}, \pier{we have that}
\begin{align}
\label{stabu4}
I_3:&=\,-\txinto\mu_2\,(g'(\rho_1)-g'(\rho_2))\,\pt\rho_1\,\pt\mu\dx\ds\,\le\,
C\tint\|\pt\rho_1(s)\|_4\,\|\rho(s)\|_4\,\|\pt\mu(s)\|_2\ds\non\\
&\le\,\frac 16\txinto|\pt\mu|^2\dx\ds\,+\,C\,\max_{0\le s\le t}\,\|\rho(s)\|_V^2
\tint\|\pt\rho_1(s)\|_V^2\ds\non\\
&\le\,\frac 16\txinto|\pt\mu|^2\dx\ds\,+\,C\tginto|\ug|^2\dg\ds\,.
\end{align}
Moreover, thanks to \eqref{stabu1} and Young's inequality, \pier{we see that}
\begin{align}\label{stabu5}
I_4:=&-\txinto \mu_2\,g'(\rho_2)\,\pt\rho\,\pt\mu\dx\ds\,\le\,\frac 16\txinto
|\pt\mu|^2\dx\ds\,+\, C\txinto|\pt\rho|^2\dx\ds\non\\
&\le \frac 16\txinto |\pt\mu|^2\dx\ds\,+\, C\,\|\ug\|_{L^2(0,t;\Hg)}^2\,.
\end{align}
Finally, we use \eqref{addregmu}, \eqref{diffi1}, \eqref{stabu1}, and H\"older's \pier{inequality to} conclude that
\begin{align}
\label{stabu6}
I_2:&=\,-\,2\txinto(g(\rho_1)-g(\rho_2))\,\pt\mu_2\,\pt\mu\dx\ds\,\le\,C\tint\|\pt\mu_2(s)\|_6
\,\|\pt\mu(s)\|_2\,\|\rho(s)\|_3\ds\non\\
&\le\,\frac 16\txinto|\pt\mu|^2\dx\ds\,+\,C\tint\|\pt\mu_2(s)\|_6^2\,\|\rho(s)\|^2_V\ds\non\\
&\le \frac 16\txinto|\pt\mu|^2\dx\ds\,+\,C\,\max_{0\le s\le t}\,\|\rho(s)\|^2_V
\tint\|\pt\mu_2(s)\|_6^2\,\ds\non\\
&\le \frac 16\txinto|\pt\mu|^2\dx\ds\,+\,C\,\|\ug\|^2_{L^2(0,t;\Hg)}\,.
\end{align}
Thus, \juerg{combining \eqref{stabu2} with \eqref{stabu3}--\eqref{stabu6}},
we have shown the estimate
\begin{align}\label{stabu7}
&\frac 13\txinto|\pt\mu|^2\dx\ds\,+\,\frac 12\,\|\mu(t)\|_V^2\non\\[1mm]
&\le\, C\,\|\ug\|^2_{L^2(0,t;L^2(\Gamma))}\,
+\,C\tint\juerg{\bigl(1+\|\pt\rho_1(s)\|_V^2\bigr)}\,\|\mu(s)\|_V^2\ds\,,
\end{align}
where the mapping $\,\,s\mapsto \|\pt\rho_1(s)\|_V^2\,\,$ is known
to belong to $L^1(0,T)$. We may therefore employ Gronwall's lemma to infer that
\beq
\label{stabu8}
\|\mu\|_{H^1(0,t;H)\cap L^\infty(0,t;V)}\,\le\,C\,\|\ug\|_{L^2(0,t;\Hg)}\,.
\eeq
It then easily follows by comparison in \eqref{diff1} that also
\begin{align}
\label{stabu9}
\|\Delta\mu\|_{\lzht}
&\gianni{{}\leq C\bigl(
  \|\pt\mu\|_{L^2(0,t;H)}
  + \|\mu\|_{L^\infty(0,t;V)}
  + \|\rho\|_{L^\infty(0,t;V)}
  + \|\pt\rho\|_{L^2(0,t;H)}
\bigr)}{} \non \\
&\leq C\,\|\ug\|_{L^2(0,t;\Hg)}\,,
\end{align}
whence, by virtue of standard elliptic estimates,
\beq\label{stabu10}
\|\mu\|_{L^2(0,t;W)}\,\le\,C\,\|\ug\|_{L^2(0,t;\Hg)}\,.
\eeq
This concludes the proof of the assertion.\qed

\section{Fr\'echet differentiability\\ of the control-to-state operator}
\setcounter{equation}{0}

\noindent In this section, we establish a differentiability result for the 
control-to-state operator $\sol$. To this end, we fix some $\bug\in {\cal U}_R$ and
set $(\bm,\br,\brg)=\sol(\bug)$, which implies that $(\bm,\br,\brg)$
satisfies \eqref{ssbounds1}, \eqref{ssbounds2}, \eqref{addregmu}, and 
$\brg=\br\onSigma$ \pier{{}a.\,e.\ on $\Sigma$}. We then consider for a fixed perturbation
$h\in{\cal X}$ \gianni{(see~\eqref{defX})} the linearized system
\begin{align}
\label{ls1}
&(1+2g(\br))\,\pt\eta\,+\,g'(\br)\pt\br\,\eta\,-\,\Delta\eta\non\\[1mm]
&\quad =\,-\,2g'(\br)\,\pt\bm\,\zeta\,-\,\bm\,g''(\br)\,\pt\br\,\zeta\,-\,\bm\,g'(\br)\,\pt\zeta
\quad\mbox{a.\,e. in }\,Q,\\[1mm]
\label{ls2}
&\pn\eta\,=\,0\quad\mbox{a.\,e. on }\,\Sigma,\quad \eta(0)=0\quad\aeO,\\[1mm]
\label{ls3}
&\pt\zeta-\Delta\zeta + \left(f''(\br)+\pi'(\br)-\bm\,g''(\br)\right)\zeta\,=\,
g'(\br)\,\eta\quad\mbox{a.\,e. in }\,Q,\\[1mm]
\label{ls4}
&\pn\zeta+\pt\zetg-\delg\zetg+\left(\fgasz(\brg)+\pigs(\brg)\right)\zetg=h,
\quad\zetg=\zeta\onSigma,
\quad\mbox{a.\,e. on }\,\Sigma,\\[1mm]
\label{ls5}
&\zeta(0)=0\quad\mbox{a.\,e. in }\,\oma,\quad\zetg(0)=0\quad
\mbox{a.\,e. on }\,\Gamma.
\end{align}
Provided that the system \eqref{ls1}--\eqref{ls5} has for every $\pier{h}\in{\cal X}$
a unique solution triple $(\eta,\zeta,\zetg)$, 
we expect that the Fr\'echet derivative $\,D\sol(\bug)\,$ of $\,\sol\,$
at $\,\bug\,$ (if it exists)
ought to be given by $\, D\sol(\bug)(h)=(\eta,\zeta,\zetg)$.  
In the following existence and uniqueness result, we show that the linearized problem
is even solvable if only $h\in L^2(\Sigma)$.

\vspace{7mm}\noindent
{\sc Theorem 3.1:} \quad{\em Suppose that} (A1)--(A6) {\em are satisfied. Then the system} \eqref{ls1}--\eqref{ls5}
{\em has for every $h\in L^2(\Sigma)$ a unique solution $(\eta,\zeta,\zeta_\Gamma)$ such that}
\begin{align}
\label{regeta}
&\eta\in H^1(0,T;H)\cap C^0([0,T];V)\cap L^2(0,T;W),\\
\label{regzeta}
&\zeta\in H^1(0,T;H)\cap C^0([0,T];V)\cap L^2(0,T;\Hdue),\\
\label{regzetag}
&\zeta_\Gamma\in H^1(0,T;\Hg)\cap C^0([0,T];\Vg)\cap L^2(0,T;H^2(\Gamma)).
\end{align}
{\em Moreover, the linear mapping $\,h\mapsto (\eta,\zeta,\zeta_\Gamma)\,$ is continuous as a mapping
from $\,L^2(\Sigma)\,$ into the Banach space}
\begin{align*}
{\cal Z}:=&\left\{(\mu,\rho,\rg)\in (H^1(0,T;H)\cap L^2(0,T;W))\times (H^1(0,T;H)\cap L^2(0,T;\Hdue))
\right.\\
&\left.\quad {}\times (H^1(0,T;\Hg)\cap L^2(0,T;H^2(\Gamma))): \,\rho\onSigma=\rg \,\,\mbox{ a.\,e. on }\,\Sigma\right\}\,.
\end{align*}

\vspace{3mm}\noindent
{\sc Proof:} \quad\, We use an approximation scheme based on a retarded argument method.
To this end, we define for every $\tau\in (0,T)$ the translation operator $\,\TT:C^0([0,T];H)
\to C^0([0,T];H)$ by setting, for all $v\in C^0([0,T];H)$,
\beq
\label{defTT}
\TT(v)(t):=v(t-\tau) \quad\mbox{if $\,t>\tau$}\quad\mbox{and }\,\TT(v)(t):=v(0) \quad
\mbox{if }\,t\le\tau.
\eeq
Notice that for every $v\in H^1(0,T;H)$ it holds that
\begin{align}
\label{TT1}
\|\TT(v)\|_{L^2(Q_t)}^2\,&\le\,\left\{
\begin{array}{l}
\|v\|_{L^2(Q_t)}^2\,+\,\tau\,\|v(0)\|^2_H
\quad\mbox{for all $t\in [\tau,T]$},\\[1mm]
t\,\|v(0)\|_H^2\quad\mbox{for all $t\in [0,\tau]$},
\end{array}
\right. \\[2mm]
\label{TT2}
\|\pt\TT(v)\|_{L^2(Q_t)}^2\,&\le\,\|\pt v\|_{L^2(Q_t)}^2\quad\mbox{for a.\,e. 
$\,t\in (0,T)$},
\end{align}
while for every $v\in C^0([0,T];V)$ we have 
\begin{align}
\label{TT3}
\|\nabla\TT(v)\|_{L^2(Q_t)}^2\,&\le\,\left\{
\begin{array}{l}
\|\nabla v\|_{L^2(Q_t)}^2\,+\,\tau\,\|\nabla v(0)\|^2_H
\quad\mbox{for all $t\in [\tau,T]$},\\[1mm]
t\,\|\nabla v(0)\|_H^2\quad\mbox{for all $t\in [0,\tau]$}.
\end{array}
\right. 
\end{align}

Now, let $N\in\nz$ be fixed, $\taun:=\gianni T/N$, as well as $t_n:=n\taun$
and $I_n:=(0,t_n)$, for $0\le n\le N$. We
then consider for every $\,n\in\{1,\ldots,N\}$ the initial-boundary value problem
\begin{align}
\label{lt1}
&(1+2g(\br))\,\pt\etan\,-\,\Delta\etan\,=\,-\,\TTN(\eta_{n-1})\,g'(\br)\pt\br\,
-\,2g'(\br)\,\pt\bm\,\zetan\non\\
&\hspace*{47.8mm}\,-\,\bm\,g''(\br)\,\pt\br\,\zetan\,-
\,\bm\,g'(\br)\,\pt\zetan
\quad\mbox{a.\,e. in }\,\oma\times I_n,\\[1mm]
\label{lt2}
&\pn\etan\,=\,0\quad\mbox{a.\,e. on }\,\Gamma\times I_n,
\quad \etan(0)=0\quad\aeO,\\[1mm]
\label{lt3}
&\pt\zetan-\Delta\zetan + \left(f''(\br)+\pi'(\br)-\bm\,g''(\br)\right)\zetan\,=\,
g'(\br)\,\TTN(\eta_{n-1})\quad\mbox{a.\,e. in }\,\oma\times I_n,\\[1mm]
\label{lt4}
&\pn\zetan+\pt\zetang-\delg\zetang+\left(\fgasz(\brg)+\pigs(\brg)\right)\zetang=h,
\quad\zetang=\zeta_{n_{|\Gamma\times I_n}},
\quad\mbox{a.\,e. on }\,\Gamma\times I_n,\\[1mm]
\label{lt5}
&\zetan(0)=0\quad\mbox{a.\,e. in }\,\oma,\quad\zetang(0)=0\quad
\mbox{a.\,e. on }\,\Gamma.
\end{align}
Here, we notice that the operator $\,\TTN\,$ acts on functions that are not defined
on the whole of $\oma\times(0,T)$; however, its meaning is still given by \eqref{defTT} if 
$n>1$, while for $n=1$ we simply set $\,\TTN(\eta_{n-1})=0$. 

The plan of the upcoming proof is as follows: in the first step, we show
that the above initial-boundary value problems have unique solutions $(\etan,\zetan,
\zetang)$ \pier{for $n=1, \ldots, N$ with the regularity as in} \eqref{regeta}--\eqref{regzetag}. Once this will
be shown, we can infer from the uniqueness~that
$$\eta_{N_{|\oma\times I_{N-1}}}=\eta_{N-1},\quad \zeta_{N_{|\oma\times I_{N-1}}}
=\zeta_{N-1}, $$
which then entails that, for almost every $(x,t)\in Q$, 
$$
\TTN(\eta_{N-1})(x,t)=\eta_{N-1}(x,t-\taun)=\eta_N(x, t-\taun)=\TTN(\eta_N)(x,t).
$$
It then follows that $(\eta^\tau,\zeta^\tau,\zeta_\Gamma^\tau):=
(\eta_N,\zeta_N,\zeta_{N_\Gamma})$ is for
$\tau=\taun$ the unique solution
to the retarded initial-boundary value problem
\begin{align}
\label{ltau1}
&(1+2g(\br))\,\pt\eta^\tau\,+\,\TT(\eta^\tau)\,g'(\br)\pt\br\,-\,\Delta
\eta^\tau\non\\[1mm]
&\quad =\,-\,2g'(\br)\,\pt\bm\,\zeta^\tau\,-\,\bm\,g''(\br)\,\pt\br\,\zeta^\tau\,-
\,\bm\,g'(\br)\,\pt\zeta^\tau
\quad\mbox{a.\,e. in }\,Q,\\[1mm]
\label{ltau2}
&\pn\eta^\tau\,=\,0\quad\mbox{a.\,e. on }\,\Sigma,
\quad \eta^\tau(0)=0\quad\aeO,\\[1mm]
\label{ltau3}
&\pt\zeta^\tau-\Delta\zeta^\tau + \left(f''(\br)+\pi'(\br)-\bm\,g''(\br)\right)
\zeta^\tau\,=\,
g'(\br)\,\TTN(\eta^\tau)\quad\mbox{a.\,e. in }\,Q,\\[1mm]
\label{ltau4}
&\pn\zeta^\tau+\pt\zeta_\Gamma^\tau-\delg\zeta_\Gamma^\tau+
\left(\fgasz(\brg)+\pigs(\brg)\right)
\zeta_\Gamma^\tau=h,
\quad\zeta_\Gamma^\tau=\zeta^\tau\onSigma,
\quad\mbox{a.\,e. on }\,\Sigma,\\[1mm]
\label{ltau5}
&\zeta^\tau(0)=0\quad\mbox{a.\,e. in }\,\oma,\quad\zeta^\tau_\Gamma(0)=0\quad
\mbox{a.\,e. on }\,\Gamma.
\end{align}

Once the unique solvability of \eqref{ltau1}--\eqref{ltau5} will be shown
for $\,\tau=\taun$, $N\in\nz$, 
\pier{in the second step of this proof we will} establish sufficiently strong a priori 
estimates, which are uniform with respect to $N\in\nz$, and then pass to
the limit as $N\to\infty$ by compactness arguments to \pier{show} the existence
of a solution $(\eta,\zeta,\zetg)$ having the required regularity properties.
As a byproduct of our estimates, we will obtain the uniqueness of the
solution and the continuity of the mapping $h\mapsto (\eta,\zeta,\zetg)$.

Pursuing our plan, we first establish the unique solvability of \eqref{lt1}--\eqref{lt5}
for every $n\in \{1,\ldots,N\}$. To this end, we argue by induction. Since the proof
for $n=1$ is similar to that used in the induction step $n-1\longrightarrow n$, we 
may confine ourselves to just perform the latter.

So let $1<n\le N$, and assume that for $1\le k\le n-1$ unique solutions 
$(\eta_k,\zeta_k,\zeta_{k_\Gamma})$ to the system \eqref{lt1}--\eqref{lt5} 
have already been constructed that satisfy for $1\le k\le n-1$ the conditions
\begin{align}
\label{regutau}
&\eta_k\in H^1(I_k;H)\cap C^0(\bar I_k;V)\cap L^2(I_k;W),\non\\
&\zeta_k\in H^1(I_k;H)\cap C^0(\bar I_k;V)\cap L^2(I_k;\Hdue),\non\\
&\zeta_{k_\Gamma}\in H^1(I_k;\Hg)\cap C^0(\bar I_k;\Vg)\cap L^2(I_k;H^2(\Gamma)).
\end{align}

First, we apply \gianni{Lemma~2.3} to infer that the initial-boundary value problem \eqref{lt3}--\eqref{lt5}
has a unique solution pair with $\,\zetan\in H^1(I_n;H)\cap C^0(\bar I_n;V)
\cap L^2(I_n;H^2(\oma))\,$ and
$\zetang\in H^1(I_n;\Hg)\cap C^0(\bar I_n;\Vg)\cap L^2(I_n;H^2(\Gamma))$. We then insert $\,\zetan\,$ in 
\eqref{lt1}. Obviously, we can rewrite the resulting identity in the form
\beq
\pt\etan\,-\,\mbox{$\frac 1{1+2g(\br)}$}\,\Delta\etan\,=\,z,
\eeq
where $\,1+2g(\br)\in \CQ$, and where, owing to \eqref{regmu}, \eqref{regrho},
and \eqref{addregmu}, the right-hand side $\,z\,$ is easily seen to belong to $\,L^2(\oma\times (0,I_n))$. 
It thus follows from
maximal parabolic regularity theory (see, e.\,g., \cite[Thm.~2.1]{DHP}) that the initial-boundary value
problem \eqref{lt1}--\eqref{lt2} enjoys a unique solution $\,\etan\in H^1(I_n;H)\cap C^0(\bar I_n;V)
\cap L^2(I_n;W)$. 
 
Now that the unique solvability of the retarded problem \eqref{ltau1}--\eqref{ltau5} with the requested
regularity is shown for every $\tau_N=T/N$, $N\in\nz$, we aim to derive a number of a priori estimates
that are uniform in $N\in\nz$. In this process, we denote by $C>0$
constants that may depend on the data of the state system but not on $N\in\nz$. For the sake of
a better readability, we will suppress the superscript \gianni{$\tau$ or $\tau_N$} during the estimations,
writing it only at the very end of each step. We also make repeated use of the global
estimates \eqref{ssbounds1}, \eqref{ssbounds2} and of \eqref{addregmu} without further
reference.

\vspace{3mm}\noindent
\underline{\sc First estimate:}
\par\nobreak
\vspace{1mm}\noindent
We add $\,\eta\,g'(\br)\,\pt\br\,$ to both sides of \eqref{ltau1} and observe that we have 
$\,\pt((\frac 12+g(\br))\,\eta^2)\,=\,(1+2g(\br))\,\eta\,\pt\eta\,+\,g'(\br)\,\pt\br\,\eta^2$.
Therefore, multiplying by $\,\eta\,$ and integrating over $Q_t$, where $0<t\le T$, 
and recalling that $\,g\,$ is nonnegative, we find that
\beq
\label{estn1}
\frac 12\xinto|\eta(t)|^2\dx
\,+\txinto|\nabla\eta|^2\dx\ds\,\le\,\sum_{j=1}^3\,I_j\,,
\eeq
where the expressions $I_j$, $1\le j\le 3$, will be specified end estimated below. At first, 
employing H\"older's inequality, \eqref{TT1}, \eqref{TT3}, and the continuity
of the embedding $V\subset L^4(\oma)$, we obtain from Young's inequality that
\begin{align}
\label{estn2}
I_1:&=\txinto g'(\br)\,\pt\br\left(\eta-\TTN(\eta)\right)\eta\dx\ds\non\\
&\le\,C\!\tint\! \|\pt\br(s)\|_4\left(\|\eta(s)\|_4+\|\TTN(\eta(s))\|_4\right)\|\eta(s)\|_2\ds
\non\\
&\le\,\frac 14\tint\|\eta(s)\|_V^2\ds\,+\,C\tint\|\pt\br(s)\|_V^2\,\|\eta(s)\|_H^2\ds\,.
\end{align}
By the same token, \pier{we have that}
\begin{align}
\label{estn3}
I_2:&=\,-\txinto\left(2g'(\br)\,\pt\bm\,+\,\bm\,g''(\br)\,\pt\br\right)\zeta\,\eta\dx\ds
\non\\
&\le\,C\tint\left(\|\pt\bm(s)\|_6+\|\pt\br(s)\|_6\right)\|\zeta(s)\|_2\,\|\eta(s)\|_3\ds\non\\
&\le\,\frac 14\tint\|\eta(s)\|_V^2\ds\,+\,C\tint\left(\|\pt\bm(s)\|_6^2+\|\pt\br(s)\|_V^2
\right)\|\zeta(s)\|_H^2\ds\,.
\end{align}
Moreover, \pier{from} Young's inequality \pier{it follows that}
\beq
\label{estn4}
I_3:=\,-\txinto\bm\,g'(\br)\,\pt\zeta\,\eta\dx\ds
\,\le\,\frac 14\txinto|\pt\zeta|^2\dx\ds\,+\,C\txinto|\eta|^2\dx\ds\,.
\eeq
Hence, combining the estimates \eqref{estn1}--\eqref{estn4}, we have shown that
\begin{align}
\label{estn5}
\|\eta(t)\|_H^2\,+\,\tint\|\eta(s)\|_V^2\ds\,&\le\,\frac 12\txinto|\pt\zeta|^2\dx\ds
\,+\,C\tint\left(1+\|\pt\br(s)\|_V^2\right)\|\eta(s)\|_H^2\ds\non\\
&\quad\,+\,C\tint\left(\|\pt\bm(s)\|_6^2+\|\pt\br(s)\|_V^2
\right)\|\zeta(s)\|_H^2\ds\,,
\end{align}
where the mappings $\,s\mapsto\|\pt\br(s)\|_V^2\,$ and $\,s\mapsto\|\pt\bm(s)\|_6^2\,$ are
known to belong to $L^1(0,T)$.

Next, we observe that Lemma 2.3 can be applied to the system \eqref{ltau3}--\eqref{ltau5}, with
$\,a:=f''(\br)+\pi'(\br)-\bm\,g''(\br)\in L^\infty(Q)$, $a_\Gamma:=f''_\Gamma(\brg)
+\pigs(\brg)\in L^\infty(\Sigma)$, $\,\sigma:=g'(\br)\,\TTN(\eta)$, and $\,\sigma_\Gamma:=h$.
We then obtain from \eqref{CoSp4} the estimate
\begin{align}
\label{estn6}
&\|\zeta\|_{H^1(0,t;H)\cap C^0([0,t];V)\cap L^2(0,t;H^2(\oma))}^2
\,+\,\|\zeta_\Gamma\|_{H^1(0,t;\Hg)\cap C^0([0,t];\Vg)\cap L^2(0,t;H^2(\Gamma))}^2\non\\
&\le\,C\tint\|\eta\pier{(s)}\|^2_H\ds\,+\,C\tginto|h|^2\dg\ds\,.
\end{align}
Combining this with \eqref{estn5}, and invoking Gronwall's lemma, we have thus shown that, 
for every $t\in (0,T]$ and $N\in\nz$, 
\begin{align}
\label{estn7}
&\|\eta^{\tau_N}\|^2_{L^\infty(0,t;H)\cap\lzvt}\,+\,
\|\zeta^{\tau_N}\|_{H^1(0,t;H)\cap C^0([0,t];V)\cap L^2(0,t;H^2(\oma))}^2\non\\
&+\,\|\zeta_\Gamma^{\tau_N}\|_{H^1(0,t;\Hg)\cap C^0([0,t];\Vg)\cap L^2(0,t;H^2(\Gamma))}^2
\,\le\,C\,\|h\|^2_{L^2(0,t;\Hg)}\,.
\end{align}

\vspace{4mm}
\noindent
\underline{\sc Second estimate:}
\par\nobreak
\vspace{1mm}\noindent
We now multiply \eqref{ltau1} by $\,\pt\eta\,$ and integrate over $Q_t$, where $0<t\le T$. Since
$\,g\,$ is nonnegative, we obtain 
\beq
\label{estn8}
\txinto|\pt\eta|^2\dx\ds\,+\,\frac 12 \,\|\nabla\eta(t)\|_H^2\,\le\,\sum_{j=1}^4\,J_j,
\eeq
where the expressions $\,J_j$, $1\le j\le 4$, will be specified and estimated below. At first,
we invoke H\"older's and Young's inequalities to obtain that
\begin{align}\label{estn9}
J_1:&=\,-\txinto g'(\br)\,\pt\br \,\TTN(\eta)\,\pt\eta\dx\ds\,\le\,C\tint\|\pt\br(s)\|_4
\,\|\TTN(\eta(s))\|_4\,\|\pt\eta(s)\|_2\dx\ds\non\\
&\le\,\frac 15\txinto|\pt\eta|^2\dx\ds\,+\,C\tint\|\pt\br(s)\|_V^2\,\|\TTN(\eta(s))\|_V^2\ds\,,
\end{align} 
where the second integral on  the \rhs, which we denote by $I(t)$, can be estimated as follows: by the definition of $\TTN$, and since $\eta(0)=0$, we obviously have that $\,I(t)=0\,$ if $\,0\le t
\le\tau_N$, while for $\,\tau_N<t\le T$ it holds that
\begin{align}
\label{estn10}
  I(t)
  & = \int_{\tau_N}^t\|\pt\br(s)\|_V^2\,\|\eta(s-\tau_N)\|_V^2\ds\,=\,
  \int_0^{t-\tau_N}\|\pt\br(s+\tau_N)\|_V^2\,\|\eta(s)\|_V^2\ds\,.
\end{align}
Hence, \pier{it is clear that}
\gianni{%
\beq
  \nonumber
  I(t) = \tint \varphi(s,t)\,\|\eta(s)\|_V^2\ds
  \quad \hbox{for every $t\in[0,T]$}
\eeq
where the function  $\varphi:[0,T]^2\to\rz$
is defined (almost everywhere with respect to~$s$)~by
$$
  \varphi(s,t) := \left\{
  \begin{array}{ll}
  0 & \hbox{if $t\leq\tau_N $ and $s\in[0,T]$}\\[1mm]
  \|\pt\br(s+\tau_N)\|_V^2 & \mbox{if $t>\tau_N $ and $0\le s\le t-\tau_N $}\\[1mm]
  0 & \mbox{if  $t>\tau_N $ and $t-\tau_N< s\le T$.}
  \end{array}
  \right.
$$
On the other hand, it holds that $\varphi(s,t)\leq\overline\varphi(s)$ for every $(s,t)\in[0,T]^2$
where  
$$
\overline\varphi(s):=\left\{
\begin{array}{ll}
\|\pt\br(s+\tau_N)\|_V^2&\mbox{if }\,0\le s\le T-\tau_N\\[1mm]
0&\mbox{if }\,T-\tau_N< s\le T .
\end{array}
\right.
$$
Thus, we also have
\beq\label{estn11}
I(t)\,\le\,\tint \overline\varphi(s)\,\|\eta(s)\|_V^2\ds
  \quad \hbox{for every $t\in[0,T]$}
\eeq
and $\overline\varphi$ is obviously bounded in $L^1(0,T)$, uniformly in $N\in\nz$.}

Next, owing to H\"older's and Young's inequalities, and invoking \eqref{estn7}, we find that
\begin{align}
\label{estn12}
J_2:&=\,-2\txinto g'(\br)\,\pt\bm\,\zeta\,\pt\eta\dx\ds\,\le\,
C\tint\|\pt\bm(s)\|_6\,\|\zeta(s)\|_3\,\|\pt\eta(s)\|_2\ds\non\\
&\le\,\frac 15\txinto|\pt\eta|^2\dx\ds\,+\,C\,\max_{0\le s\le t}\,\|\zeta(s)\|^2_V
\tint\|\pt\bm(s)\|_6^2\ds\non\\
&\le\,\frac 15\txinto|\pt\eta|^2\dx\ds\,+\,C\,\|h\|^2_{L^2(0,t;\Hg)}\,,
\end{align}
as well as
\begin{align}
\label{estn13}
J_3:&=\,-\txinto \bm\,g''(\br)\,\pt\br \,\zeta\,\pt\eta\dx\ds\,\le\,
C\tint\|\pt\br(s)\|_6\,\|\zeta(s)\|_3\,\|\pt\eta(s)\|_2\ds\non\\
&\le\,\frac 15\txinto|\pt\eta|^2\dx\ds\,+\,C\,\|h\|^2_{L^2(0,t;\Hg)}\,.
\end{align}
Finally, \gianni{owing to \eqref{estn7} once more,} we obtain that
\begin{align}\label{estn14}
J_4:&=\,-\txinto\bm\,g'(\br)\,\pt\zeta\,\pt\eta\dx\ds\,\le\,\frac 15\txinto|\pt\eta|^2\dx\ds\,+\,
C\txinto|\pt\zeta|^2\dx\ds\non\\
&\le\,\frac 15\txinto|\pt\eta|^2\dx\ds\,+\,C\,\|h\|^2_{L^2(0,t;\Hg)}\,.
\end{align}
Combining the estimates \eqref{estn8}--\eqref{estn14}, we can infer from Gronwall's lemma that
\beq\label{estn15}
\left\|\eta^{\tau_N}\right\|^2_{H^1(0,t;H)\cap L^\infty(0,t;V)}\,\le\,C\,\|h\|^2_{L^2(0,t;\Hg)}\,.
\eeq
Then, \gianni{by comparing in \eqref{ltau1} and using the full regularity of~$(\bm,\br)$
(in~particular \eqref{addregmu})}, we easily check that also
\beq\label{estn16}
\left\|\Delta\eta^{\tau_N}\right\|^2_{L^2(0,t;H)}\,\le\,C\,\|h\|^2_{L^2(0,t;\Hg)}\,,
\eeq
whence, by standard elliptic estimates,
\beq\label{estn17}
\left\|\eta^{\tau_N}\right\|^2_{L^2(0,t;W)}\,\le\,C\,\|h\|^2_{L^2(0,t;\Hg)}\,.
\eeq
In conclusion, by virtue of \eqref{estn7}, \eqref{estn15}, \eqref{estn17}, and since
the embedding $\,(H^1(0,t;H)\cap L^2(0,t;\Hdue))\subset C^0([0,t];V)\,$ is continuous,
we have shown the estimate
\begin{align}
\label{estn18}
&\|\eta^{\tau_N}\|^2_{H^1(0,t;H)\cap C^0([0,t];V)\cap L^2(0,t;\Hdue)}\,+\,
\|\zeta^{\tau_N}\|_{H^1(0,t;H)\cap C^0([0,t];V)\cap L^2(0,t;H^2(\oma))}^2\non\\
&+\,\|\zeta_\Gamma^{\tau_N}\|_{H^1(0,t;\Hg)\cap C^0([0,t];\Vg)\cap L^2(0,t;H^2(\Gamma))}^2
\,\le\,C\,\|h\|^2_{L^2(0,t;\Hg)}\non\\
&\quad\mbox{for all $\,N\in\nz\,$ and $\,t\in (0,T]$} .
\end{align}

We are now in a position to show the existence of a solution to \eqref{ls1}--\eqref{ls5}.
Indeed, thanks to \eqref{estn18}, there are functions $(\eta,\zeta,\zetg)$, such that,
for a subsequence which is again indexed by $N$, we have for $N\to\infty$ that
\begin{align}
\label{coneta}
\eta^{\tau_N}&\to \eta\quad\,\,\mbox{weakly in }\,H^1(0,T;H)\cap C^0([0,T];V)
\cap L^2(0,T;W)\,,\\[1mm]
\label{conzeta}
\zeta^{\tau_N}&\to \zeta\quad\,\,\mbox{weakly in }\,H^1(0,T;H)\cap C^0([0,T];V)
\cap L^2(0,T;\Hdue)\,,\\[1mm]
\label{conzetg}
\zeta_\Gamma^{\tau_N}&\to \zetg\quad\mbox{weakly in }\,H^1(0,T;\Hg)\cap C^0([0,T];\Vg)
\cap L^2(0,T;H^2(\Gamma))\,.
\end{align} 
This implies, in particular, that the initial and boundary conditions \eqref{ls2} and 
\eqref{ls5} are fulfilled, and, since 
$\,\zeta^{\tau_N}\onSigma\to \zeta\onSigma\,$
weakly in $\,L^2(0,T;H^{3/2}(\Gamma))\,$ 
by \eqref{conzeta} and the trace theorem,
we have that $\,\zeta\onSigma=\zetg\,$ almost everywhere on $\Sigma$. 

Moreover, thanks to \cite[Sect.~8, Cor.~4]{Simon}, we may
without loss of generality assume that, for every $p\in [1,6)$,
\begin{align}
\label{conny1}
\eta^{\tau_N}&\to \eta\quad\,\,\,\mbox{strongly in } \,C^0([0,T];L^p(\oma))\,,\\[1mm]
\label{conny2}
\zeta^{\tau_N}&\to \zeta\quad\,\,\mbox{strongly in } \,C^0([0,T];L^p(\oma))\,,\\[1mm]
\label{conny3}
\zeta_\Gamma^{\tau_N}&\to \zetg\quad\mbox{strongly in } \,C^0([0,T];L^p(\Gamma))\,.
\end{align}
In addition, it holds $\,\TTN(\eta^{\tau_N})\to\eta\,$ strongly in $L^2(Q)$, and it is
easily verified that
\begin{align}
\label{conny4}
&\TTN(\eta^{\tau_N})\,g'(\br)\,\pt\br \to \eta\,g'(\br)\,\pt\br, \enskip
g'(\br)\,\pt\bm\,\zeta^{\tau_N}\to g'(\br)\,\pt\bm\,\zeta,\non\\[1mm]
&\bm\,g''(\br)\,\pt\br\,\zeta^{\tau_N} \to \bm\,g''(\br)\,\pt\br\,\zeta,
\enskip \pier{g'(\br)\, \TTN(\eta^{\tau_N}) \to g'(\br)\,\eta}\, \juerg , 
\quad\mbox{all weakly in }\,L^1(Q)\,.
\end{align}
Therefore, we may pass to the limit as $N\to\infty$ in \eqref{ltau1}--\eqref{ltau5},
written for $\tau=\tau_N$, to conclude that the triple $(\eta,\zeta,\zetg)$ is in
fact a solution to the system \eqref{ls1}--\eqref{ls5} that enjoys the regularity
properties \eqref{regeta}--\eqref{regzetag}. Moreover, passage to the limit as 
$N\to\infty$ in \eqref{estn18}, using the weak sequential semicontinuity of norms, yields
that 
\begin{align}
\label{estn19}
&\|\eta\|^2_{H^1(0,t;H)\cap C^0([0,t];V)\cap L^2(0,t;\Hdue)}\,+\,
\|\zeta\|_{H^1(0,t;H)\cap C^0([0,t];V)\cap L^2(0,t;H^2(\oma))}^2\non\\
&+\,\|\zeta_\Gamma\|_{H^1(0,t;\Hg)\cap C^0([0,t];\Vg)\cap L^2(0,t;H^2(\Gamma))}^2
\,\le\,C\,\|h\|^2_{L^2(0,t;\Hg)}\,\quad\mbox{for all $\,t\in (0,T]$}\,.
\end{align}

It remains to show that the solution is unique, which, in view of \eqref{estn19},
would entail that the linear mapping $\,h\mapsto (\eta,\zeta,\zetg)\,$ is continuous from
$L^2(\Sigma)$ into ${\cal Z}$. So let us assume that two solutions $(\eta_i,
\zeta_i, \zeta_{i_\Gamma})$, $i=1,2$, satisfying \eqref{regeta}--\eqref{regzetag}
are given. Then the triple $\,(\eta,\zeta,\zetg)$, where
 $\eta:=\eta_1-\eta_2$, $\zeta:=\zeta_1-\zeta_2$, 
$\zetg:=\zeta_{1_\Gamma}-\zeta_{2_\Gamma}$, satisfies Eqs.\ \eqref{ls1}--\eqref{ls5}
with $h=0$. 

At this point, we can repeat the estimations performed in the
{\sc First estimate} above, where the only difference (which even simplifies
the analysis) is given by the fact that in Eq.\ \eqref{ls1} the term
$\,\eta\,g'(\br)\,\pt\br\,$ appears in place of the expression
$\,\TT(\eta^\tau)\,g'(\br)\,\pt\br\,$ occurring in Eq.\ \eqref{ltau1}. We thus
can claim that the estimate \eqref{estn7} is valid with
$\left(\eta^{\tau_N},\zeta^{\tau_N},\zeta_\Gamma^{\tau_N}\right)$ replaced
by $(\eta,\zeta,\zetg)$. Since $h=0$ in the present situation, we obtain that
$\eta=\zeta=0$ almost everywhere in $Q$, and $\zetg=0$ almost everywhere on
$\Sigma$. This concludes the proof of the assertion.\qed 

\vspace{5mm}
We are now in a position to prove the Fr\'echet differentiability of the control-to-state operator. 
We \gianni{recall the definition \eqref{defX} of $\cal X$ and state} the following result.

\vspace{7mm}\noindent
{\sc Theorem 3.2:} \quad\,{\em Suppose that the conditions} (A1)--(A5) {\em are fulfilled. Then the
control-to-state operator $\,\sol:\ug\mapsto (\mu,\rho,\rg)\,$ is Fr\'echet differentiable as a mapping from $\,{\cal U}_R
\subset {\cal X}\,$ into the Banach space}
\begin{align*}
{\cal Y}:=&\left\{(\mu,\rho,\rg)\in (L^\infty(0,T;H)\cap L^2(0,T;V))
\times (H^1(0,T;H)\cap L^2(0,T;\Hdue))
\right.\\
&\left.\quad{}\times (H^1(0,T;\Hg)\cap L^2(0,T;H^2(\Gamma))): \,\gianni{\rg=\rho\onSigma} \,\,\mbox{ a.\,e. on }\,\Sigma\right\}\,.
\end{align*} 
{\em Moreover, for every $\bug\in{\cal U}_R$, the Fr\'echet derivative $\,D\sol(\bug)\in {\cal L}({\cal X},
{\cal Y})\,$ is evaluated at any $h\in{\cal X}$ by 
putting $\,D\sol(\bug)(h):=(\eta,\zeta,\zeta_\Gamma)$, where $\,(\eta,\zeta,\zeta_\Gamma)\,$ is the
unique solution to the linearized system} \eqref{ls1}--\eqref{ls5}.

\vspace{3mm}\noindent
{\sc Proof:} \quad\, According to Theorem 3.1, the linear mapping $\,h\mapsto 
(\etah,\zeh,\zeg):=(\eta,\zeta,\zetg)\,$ is continuous from $L^2(\Sigma)$ into
${\cal Z}$ and thus, a fortiori, also from ${\cal X}$ into~${\cal Y}$. Hence, if
the derivative $\,D\sol(\bug)\,$ exists and has the asserted form, then it belongs 
to $\,{\cal L}({\cal X},{\cal Y})$.

Now notice that ${\cal U}_R$ is open in ${\cal X}$, and thus there is some
$\Lambda>0$ such that $\bug+h\in {\cal U}_R$ whenever $\|h\|_{{\cal X}}\le\Lambda$.
In the following, we consider only such perturbations~$h$. We then put, for 
any such~$h$, 
$$(\muh,\rh,\rgh):=\sol(\bug+h), \quad z^h:=\muh-\bm-\etah,
\quad y^h:=\rh-\br-\zeh, \quad y^h_\Gamma:=\rgh-\brg-\zeg,$$
where \pier{$(\bm,\br,\brg):=\sol(\bug)$ and} \gianni{$(\etah, \zeh,\zeg)$ denotes the unique solution $(\eta,\zeta,\zetg)$} to the linearized system \eqref{ls1}--\eqref{ls5}. Notice 
that we have $\,y^h_\Gamma=y^h\onSigma$, as well as   
\begin{align}
\label{regzh}
&z^h\in H^1(0,T;H)\cap C^0([0,T];V)\cap L^2(0,T;W),\\
\label{regyh}
&y^h\in H^1(0,T;H)\cap C^0([0,T];V)\cap L^2(0,T;\Hdue),\\
\label{regygh}
&\ygh\in H^1(0,T;\Hg)\cap C^0([0,T];\Vg)\cap L^2(0,T;H^2(\Gamma)).
\end{align}
We also notice that the global bounds \eqref{ssbounds1} and \eqref{ssbounds2}
are satisfied for $(\muh,\rh,\rgh)$, and, owing to Theorem 2.4, we have the
global stability  estimate
\begin{align}
\label{stabil2}
&\|\muh-\bm\|_{H^1(0,t;H)\cap C^0([0,t];V)\cap L^2(0,t;W)}
\,+\,\|\rh-\br\|_{H^1(0,t;H)\cap C^0([0,t];V)\cap L^2(0,t;\Hdue)}\non\\
&+\,\|\rgh-\brg\|_{H^1(0,t;\Hg)\cap C^0([0,t];\Vg)\cap L^2(0,t;H^2(\Gamma))}
\,\le\, K^*_3\,\|h\|_{L^2(0,t;\Hg)} \quad\forall\,t\in(0,T).
\end{align}
Moreover, by Taylor's theorem and \eqref{ssbounds2}, it holds that
\begin{align}
\label{taylor1}
&\bigl|f'(\rh)-f'(\br)-f''(\br)\,\zeta^h\bigr|\,+\,\bigl|g(\rh)-g(\br)-g'(\br)\,\zeta^h\bigr|\,+\,
\bigl|g'(\rh)-g'(\br)-g''(\br)\,\zeta^h\bigr|\non\\
&+\,\bigl|\pi(\rh)-\pi(\br)-\pi'(\br)\,\zeta^h\bigr|\,\le\,C\left(|y^h|\,+\,|\rh-\br|^2\right)
\quad\aeQ,\\[1mm]
\label{taylor2}
&\bigl|\fgas(\rgh)-\fgas(\brg)-\fgasz(\brg)\,\zeg\bigr|\,+\,
\bigl|\pig(\rh_\Gamma)-\pig(\br_\Gamma)-\pigs(\br_\Gamma)\,\zeta^h_\Gamma\bigr|\non\\
&\le\,C\left(|y^h_\Gamma|\,+\,|\rho^h_\Gamma
-\br_\Gamma|^2
\right)\quad\mbox{a.\,e. on }\,\Sigma,
\end{align}
where, here and in the remainder of the proof, we denote by $C>0$ constants that may 
depend on the data of the system but not on
the special choice of $h$ with $\|h\|_{\cal X}\le\Lambda$. The actual value of $\,C\,$
may change between lines and  even within formulas. 

According to the definition of the notion of Fr\'echet differentiability, we need
to show~that
\beq\label{Frechet1}
\lim_{\|h\|_{\cal X}\to 0}\,\frac{\left\|\sol(\bug+h)-\sol(\bug)-(\eta^h,
\zeh,\zeg)\right\|_{\cal Y}}{\|h\|_{\cal X}}\,=\,0\,.
\eeq
It thus suffices to \gianni{prove} the existence of an increasing function 
$\,Z:(0,\Lambda)\to (0,+\infty)\,$ such that $\,\lim_{\lambda\searrow0}\,\frac
{Z(\lambda)}{\lambda^2}=0\,$ and 
\begin{align}
\label{Frechet2}
&\|z^h\|^2_{L^\infty(0,T;H)\cap L^2(0,T;V)}\,+\,\|y^h\|^2_{H^1(0,T;H)\cap L^2(0,T;\Hdue)}
\,+\,\|\ygh\|^2_{H^1(0,T;\Hg)\cap L^2(0,T;H^2(\Gamma))}\non\\
&\le\,
Z\bigl(\|h\|_{H^1(0,T;\Hg)}\bigr)\,.
\end{align}
To begin with, using the state system \eqref{ss1}--\eqref{ss5} and the linearized
system \eqref{ls1}--\eqref{ls5}, we easily verify that the triple 
$(z^h,y^h,\ygh)$ is a strong solution to the system
\begin{align}
\label{yz1}
&(1+2g(\br))\,z^h_t\,+\,g'(\br)\br_t\,z^h\,+\,\bm\,g'(\br)\,y^h_t-\Delta \gianni z^h\nonumber\\
&=-\,2\left(g(\rh)-g(\br)\right)\left(\mu^h_t-\bm_t\right)\,
-\,2\,\bm_t\left(g(\rh)-g(\br)-
g'(\br)\zeta^h\right)\nonumber\\
&\quad \,-\,\bm\,\br_t\left(g'(\rh)-g'(\br)- g''(\br)\zeh\right)
\,-\,\bm\left(g'(\rh)-g'(\br)\right)\left(\rh_t-\br_t\right)\nonumber\\
&\quad\,-\left(\muh-\bm\right)\left[\bigl(g'(\rh)-g'(\br)\bigr)\,\br_t
+ g'(\rh)\left(\rh_t-\br_t\right)\right] \quad\mbox{a.\,e. in }\, Q,\\[1mm]
\label{yz2}
&\hspace*{10mm}\pn z^h=0 \quad\mbox{a.\,e. on }\,\Sigma,
\quad z^h(0)=0 \quad\mbox{a.\,e. in }\,\Omega,\\[1mm]
\noalign{\allowbreak}
\label{yz3}
&y^h_t-\Delta y^h=-\left(f'(\rh)-f'(\br)-f''(\br)\zeh\right)-
\left(\pi(\rh)-\pi(\br)-\pi'(\br)\,\zeh\right)
\nonumber\\
&\hspace*{22mm} +\,g'(\br)\,z^h\,+\,\bm\left(g'(\rh)-g'(\br)-g''(\br)\zeh\right)\nonumber\\
&\hspace*{22mm} +\left(\muh-\bm\right)\left(g'(\rh)-g'(\br)\right)  
\quad\mbox{a.\,e. in}\,Q,\\[1mm]
\label{yz4}
&\pn y^h\,+\,\pt \ygh\,-\,\delg\ygh\,=\,-\,\bigl(\fgas(\rho^h_\Gamma)-\fgas(\brg)
-\fgasz(\brg)\,\zeg\bigr)\non\\
&\hspace*{9mm}-\,\pier{\bigl(\pig(\rho^h_\Gamma)-\pig(\brg)-\pigs(\brg)\zeg\bigr),}
\quad y^h_\Gamma=y^h\onSigma,
 \quad\mbox{a.\,e. on }\,\Sigma,\\[1mm]
\label{yz5} 
&\hspace*{1cm} y^h(0)=0 \quad\mbox{a.\,e. in }\,\Omega, \quad
y^h_\Gamma(0)=0 \quad\mbox{a.\,e. on }\,\Gamma.
\end{align}

In the following, we make repeated use of the mean value theorem and 
of the global estimates \eqref{ssbounds1}, \eqref{ssbounds2}, and \eqref{stabil2},
without further reference. 
For the sake \revis{of better readability}, we will
omit the superscript $h$ of the quantities $z^h, y^h, y^h_\Gamma $  during the estimations, writing it only
at the end of the respective estimates.

\vspace{3mm}\noindent
\underline{\sc First estimate:}
\par\nobreak
\vspace{1mm}\noindent
Let \gianni{an arbitrary} $t\in (0,T]$ be fixed. 
First, \gianni{let us} observe that
$\,\,\pt\left((\frac 12 + g(\br))z^2\right)\,=
\,(1+2\,g(\br))\,z\,z_t \,+\,g'(\br)\,\br_t\,z^2\,.
$ 
Hence, \gianni{adding the same term $z$ to both sides of \eqref{yz1} for convenience,} 
multiplication by $z$ and  integration over $Q_t$ \gianni{yield} the estimate
\gianni{%
\beq
\label{p321}
\xinto\left(\mbox{$\frac 12$} + g(\br(t))\right)z^2(t)\dx
  + \int_0^t \|z(s)\|_V^2\ds
  \leq \int_0^t \|z(s)\|_H^2\ds
  + C\sum_{j=1}^7 |I_j|\,,
\eeq
}%
where the quantities $I_j$, $1\le j\le 7$, are specified and estimated as follows: at first, Young's
inequality shows that, for every $\gamma>0$ (to be chosen later),
\beq
\label{p322}
I_1:=\,-\txinto \bm\,g'(\br)\,y_t\,z\dx\ds\,\le\,\gamma\txinto y_t^2\dx\ds\,+\,\frac C\gamma\txinto z^2\dx\ds\,.
\eeq  
Moreover, we have, by H\"older's and Young's inequalities and \eqref{stabil2},
\begin{align}
\label{p323}
I_2:&=\,-2\txinto\left(g(\rh)-g(\br)\right)\left(\muh_t-\bm_t\right)z\dx\ds\nonumber\\[1mm]
&\le\,C\tint\|\rh(s)-\br(s)\|_6\,\|\muh_t(s)-\bm_t(s)\|_2\,\|z(s)\|_3\ds\nonumber\\[1mm]
&\le\,
C\,\|\rh-\br\|_{C^0([0,t];V)}\,\|\muh-\bm\|_{H^1(0,t;H)}\,\|z\|_{\lzvt}\nonumber\\[1mm]
&\le\,\gamma\,\|z\|_{\lzvt}^2\,+\,\frac C\gamma \,\|h\|^4_{L^2(0,t;\Hg)}\,.
\end{align}
Next, we employ \eqref{taylor1}, the H\"older and Young inequalities, and \eqref{stabil2}, to infer that
\begin{align}
\label{p324}
I_3:&=\,-2\txinto\bm_t\left(g(\rh)-g(\br)- g'(\br)\zeta^h\right)z\dx\ds\nonumber\\[1mm]
&\le \,C\txinto|\bm_t|\left(|y|\,+\,|\rh-\br|^2\right)\,|z|\dx\ds\nonumber\\[1mm]
&\le \,C\tint\|\bm_t(s)\|_6\left(\|y(s)\|_3\,\|z(s)\|_2\,+\,\|\rh(s)-\br(s)\|_6^2
\,\|z(s)\|_2\right)\ds\nonumber\\[1mm]
\noalign{\allowbreak}
&\le\,C\tint\|\bm_t(s)\|_6^2\,\|z(s)\|_H^2\ds
\,+\,C\tint\|y(s)\|_V^2\ds\,+\,C\tint\|\rh\pier{(s)}-\br\pier{(s)}\|_V^4\ds\nonumber\\[1mm]
&\le \,C\tint\gianni{\bigl(1+\|\bm_t(s)\|_6^2\bigr)}
\left(\|y(s)\|_V^2\,+\,\|z(s)\|^2_H\right)ds\,+\,C\,\|h\|_{L^2(0,t;\Hg)}^4\,.
\end{align}
Likewise, with \eqref{ssbounds2}, \eqref{taylor1}, \eqref{stability}, and the H\"older and Young 
inequalities, we find that
\begin{align}
\label{p325}
I_4:&=-\!\txinto\bm\,\br_t \,(g'(\rh)-g'(\br)- g''(\br)\zeta^h)\,z
\dx\ds\le C\!\txinto|\br_t|(|y|+|\rh-\br|^2)|z|\dx\ds\nonumber\\
&\le\,C\tint\|\br_t(s)\|_6\left(\|y(s)\|_3\,+\,\|\rh(s)-\br(s)\|_6^2\right)\|z(s)\|_2
\dx\ds\nonumber\\[1mm]
&\le\,C\tint\|y(s)\|_V^2\ds\,+\,C\tint\|\br_t(s)\|_V^2\,\|z(s)\|^2_H\ds
\,+\,C\,\max_{0\le s\le t}\,\|\rh(s)-\br(s)\|_V^4\nonumber\\[1mm]
&\le \,C\tint\|y(s)\|_V^2\ds\,+\,C\tint\|\br_t(s)\|_V^2\,\|z(s)\|^2_H\ds
\,+\,C\,\|h\|^4_{L^2(0,t;H_\Gamma)}\,.
\end{align}

In addition, arguing similarly, we have
\begin{align}
\label{p326}
I_5:&=\,-\txinto\bm\left(g'(\rh)-g'(\br)\right)\left(\rh_t-\br_t\right)z\dx\ds
\nonumber\\[1mm]
&\le\,C\tint\|\rh(s)-\br(s)\|_6\,\|\rh_t(s)-\br_t(s)\|_2\,\|z(s)\|_3\ds\nonumber\\[1mm]
&\le\,C\,\|\rh-\br\|_{C^0([0,t];V)}
\,\|\rh-\br\|_{H^1(0,t;H)}\,\|z\|_{\lzvt}\nonumber\\[1mm]
&\le\,\gamma\tint\|z(s)\|_V^2\ds\,+\,\frac C\gamma\,\|h\|^4_{L^2(0,t;H_\Gamma)}\,,
\end{align}
as well as
\begin{align}
\label{p327}
I_6:&=\,-\txinto\br_t\left(\muh-\bm\right)\left(g'(\rh)-g'(\br)\right)z\dx\ds\nonumber\\[1mm]
&\le\,C\tint\|\br_t(s)\|_6\left\|\muh(s)-\bm(s)\right\|_6\left\|\rh(s)-\br(s)\right\|_6
\,\|z(s)\|_2\dx\ds\nonumber\\[1mm]
&\le\,C\tint\|\br_t(s)\|_V^2\,\|z(s)\|_H^2\ds\,
 +\,C\,\gianni{\|\muh-\bm\|_{C^0([0,t];V)}^2\,\|\rh-\br\|_{C^0([0,t];V)}^2}\nonumber\\[1mm]
&\le\,C\tint\|\br_t(s)\|_V^2\,\|z(s)\|_H^2\ds\, +\,C\,\|h\|^4_{L^2(0,t;\Hg)}\,.
\end{align}
Finally, we find that
\begin{align}
\label{p328}
I_7:&=\,-\txinto \left(\muh-\bm\right)g'(\rh)\left(\rh_t-\br_t\right)z\dx\ds\nonumber\\[1mm]
&\le\,C\tint\|\muh(s)-\bm(s)\|_6\,\|\rh_t(s)-\br_t(s)\|_2\,\|z(s)\|_3\ds\nonumber\\[1mm]
&\le\,C\,\|\muh-\bm\|_{C^0([0,t];V)}\,\|\rh-\br\|_{H^1(0,t;H)}\,\|z\|_{\lzvt}\nonumber\\[1mm]
&\le\,\gamma\tint\|z(s)\|_V^2\ds\,+\,\frac C\gamma\,\|h\|^4_{L^2(0,t;H_\Gamma)}\,.
\end{align}

In conclusion, combining the estimates \eqref{p321}--\eqref{p328}, 
and choosing $\gamma=\frac 18$,
we have shown that
\begin{align}
\label{p329}
&\frac 12\,\left\|z^h(t)\right\|_H^2\,+\,\frac 12\tint\left\|z^h(s)\right\|^2_V ds\,
\le\,\frac 18\txinto\left|y^h_t\right|^2 dx\ds\,+
\,C\,\|h\|^4_{L^2(0,t;H_\Gamma)}\nonumber\\[1mm]
&\hspace*{1cm} +\,C\tint\left(1+\left\|\bm_t(s)\right\|_6^2\,+
\left\|\br_t(s)\right\|_V^2\right)
\left(\left\|y^h(s)\right\|_V^2
\,+\left\|z^h(s)\right\|_H^2\right)\ds\,,
\end{align}
where we observe that, in view of \eqref{regmu} and \eqref{addregmu}, the mapping 
$\,s\mapsto \|\bm_t(s)\|_6^2\,+\,\|\br_t(s)\|_V^2\,$ belongs to
$L^1(0,T)$. 

\vspace{3mm}\noindent
\underline{\sc Second estimate:} \,\,\,We now observe that $\,y^h\,$ satisfies a linear problem
of the form \eqref{CoSp1}--\eqref{CoSp3}, where in this case $\,\pier{a=0}\,$ and $\,\pier{a_\Gamma=0}$,
and where $\,\sigma\,$ and $\sigma_\Gamma\,$ are equal to the right-hand sides of 
\eqref{yz3} and \eqref{yz4}, respectively. We therefore have, with this choice of $\sigma,
\sigma_\Gamma$,
\begin{align}
\label{neu1}
&\|y^h\|_{H^1(0,t;H)\cap C^0([0,t];V)\cap L^2(0,t;H^2(\oma))}
\,+\,\|y^h_\Gamma\|_{H^1(0,t;\Hg)\cap C^0([0,t];\Vg)\cap L^2(0,t;H^2(\Gamma))}\non\\
&\le\,C_L\left(\|\sigma\|_{L^2(Q_t)}\,+\,\|\sigma_\Gamma\|_{L^2(\Sigma_t)}\right)\quad\forall\,t\in (0,T]\,.
\end{align}
Now, using \eqref{taylor1}, \eqref{taylor2}, and the stability estimate \eqref{stability}, we easily
conclude that
\begin{align}
\|\sigma\|_{L^2(Q_t)}^2&\le\,C\txinto\left(|y^h|^2\,+\,|z^h|^2\,+\,|\rh-\br|^4\,
+\,|\muh-\bm|^2\,|\rh-\br|^2\right)\dx\ds\non\\
&\le\,C\txinto\left(|y^h|^2+|z^h|^2\right)\dx\ds\,+\,C\,\|h\|^4_{L^2(0,t;H_\Gamma)}\,,\label{neu2}\\[1mm]
\|\sigma_\Gamma\|_{L^2(\Sigma_t)}^2&\le\,C\tginto\left|y^h_\Gamma\right|^2\dg\ds\,+\,C\tginto
\left|\rho^h_\Gamma-\br_\Gamma\right|^4\dg\ds\non\\
&\le C\tginto\left|y^h_\Gamma\right|^2\dg\ds\,+\,C\,\|h\|^4_{L^2(0,t;H_\Gamma)}\,.\label{neu2bis}
\end{align}
Thus, combining the estimates \pier{\eqref{p329}--\eqref{neu2bis}} and invoking Gronwall's lemma,
we have proved the estimate
\begin{align}
\label{neu3}
&\|z^h\|_{C^0([0,t];H)\cap L^2(0,t;V)}^2 \,+\, \|y^h\|_{H^1(0,t;H)\cap C^0([0,t];V)\cap
 L^2(0,t;H^2(\oma))}^2\non\\[1mm]
&+\,\|y^h_\Gamma\|_{H^1(0,t;\Hg)\cap C^0([0,t];\Vg)\cap L^2(0,t;H^2(\Gamma))}^2\,\le\,
\gianni{\widetilde C}\,\|h\|_{L^2(0,t;H_\Gamma)}^4
\end{align}
\gianni{where $\,\widetilde C\,$ is a sufficiently large constant.
Therefore, the condition \eqref{Frechet2} is satisfied for the function $\,Z(\lambda)=\widetilde C\,\lambda^4$.}
This concludes the proof of the assertion.\qed

\vspace{7mm}
We are now in the position to state the following necessary optimality condition, 
which is a simple standard application of
the chain rule and of the fact that $\uad$ is a convex set. We thus may leave its proof to the reader. 

\vspace{5mm}\noindent
{\sc Corollary 3.3:}  \quad{\em Let the general hypotheses} (A1)--(A6) {\em be
fulfilled, and assume that $\bug\in\uad$ is a solution to the control problem} {\bf (CP)}
{\em with associated state $(\bm,\br,\br_\Gamma)=\sol (\bug)$. Then we have, for every 
$v_\Gamma\in \uad$,}
\begin{align}
\label{vug1}
&\beta_1\texinto (\bm-\hat\mu_Q)\,\eta\dx\dt\,+\,\beta_2\texinto(\br-\hat\rho_Q)\,\zeta\dx\dt
\,+\,\beta_3\teginto(\br_\Gamma-\hat\rho_\Sigma)\,\zeta_\Gamma\dg\dt
\nonumber\\[1mm]
&+\,\beta_4\xinto(\br(T)-\hat\rho_\oma)\,\zeta(T)\dx\,+\,\beta_5\ginto(\br_\Gamma(T)-\hat\rho_\Gamma)
\,\zeta_\Gamma(T)\dg\nonumber\\[1mm]
&+\,\beta_6\teginto\bug\,(v_\Gamma-\bug)\dg\dt\,\ge\,0\,,
\end{align}
{\em where $(\eta,\zeta,\zeta_\Gamma)$ denotes the (unique) solution to the linearized system} 
\eqref{ls1}--\eqref{ls5} {\em associated with \,$h=v_\Gamma-\bug$.}

\section{Existence and necessary optimality conditions}
\setcounter{equation}{0}

In this section, we state and prove the main results of this paper. We begin with an
existence result.

\vspace{3mm}\noindent
{\sc Theorem 4.1:} \quad\,{\em Suppose that the conditions} (A1)--(A6) 
{\em are fulfilled. Then the optimal control problem} {\bf (CP)} {\em admits a solution
$\ug\in\uad$.} 

\vspace{3mm}
\noindent
{\sc Proof:} \quad\,Since $\uad\not =\emptyset$, we may pick a minimizing sequence
$\,\{u_{\Gamma,n}\}_{n\in\nz}\subset\uad\,$ for the control problem. Now put
$(\mu_n,\rho_n,\rho_{n_\Gamma}):=\sol(u_{\Gamma,n})$, where $\,\rho_{n_\Gamma}=
\rho_{n\onSigma}$, for $n\in\nz$. By virtue of the global estimates \pier{\eqref{ssbounds1},
\eqref{ssbounds2}} and of the separation property \eqref{separ}, 
and invoking \cite[Sect.~8, Cor.~4]{Simon}, we may without loss of
generality assume that there exist some $\bug\in\uad$ and functions $\,\bm,\br,\brg$ such that,
as $\,n\to\infty$,
\begin{align}
\label{con1}
u_{\Gamma,n}&\to \bug\quad\mbox{weakly-star in $H^1(0,T;\Hg)\cap L^\infty(\Sigma)$}, \\[1mm]
\label{con2}
\mu_n&\to\bm \quad\hspace*{1.5mm}
\pier{\mbox{weakly star in $ H^1(0,T;H)\cap L^\infty(0,T;V)\cap L^2(0,T;W)\cap L^\infty (Q)$}} \non \\[1mm]
&\hspace*{14.5mm}\pier{\mbox{and strongly in }\,C^0([0,T];H)\cap L^2(0,T;L^\infty (\Omega)),}\\[1mm]
\label{con3}
\rho_n&\to \br \quad\hspace*{1.8mm}\mbox{weakly-star in }\, W^{1,\infty}(0,T;H)
\cap H^1(0,T;V)\cap L^\infty(0,T;\Hdue)\non\\
&\hspace*{14.5mm} \mbox{and strongly in }\,\CQ  ,
\\[1mm]
\label{con4}
\rho_{n_\Gamma}&\to\brg\hspace*{3.5mm} \mbox{weakly-star in }\,
 W^{1,\infty}(0,T;\Hg) \cap H^1(0,T;\Vg)\cap L^\infty(0,T;H^2(\Gamma)),\\[1mm]
 &\hspace*{15mm} -1<r_*\le\rho_n(x,t)\le r^*<+1\quad\forall\,(x,t)\in\overline Q\,.
\end{align}
In particular, it holds \,\,$\rho_{n_\Gamma}=\rho_{n\onSigma}\to \br\onSigma\,$ strongly
in $\,\CS$, which entails that $\,\brg=\br\onSigma\,$ on $\overline\Sigma$ and, thanks to 
the assumptions
(A2) and (A3), that
\begin{align}
\label{con5}
\Phi(\rho_n)&\to \Phi(\br)\hspace*{8mm}\mbox{strongly in $\,\CQ$\, for }\,\Phi\in\{g,g',f',\pi\},\\[1mm]
\label{con6}
\Phi_\Gamma(\rho_{n_\Gamma})&\to \Phi_\Gamma(\brg)
\hspace*{4mm}\mbox{strongly in $\,\CS$\, for }\,
\Phi_\Gamma\in\{f'_\Gamma,\pig\}.
\end{align}
Moreover, owing to the trace theorem,
\beq\label{con7}
\pn\mu_n\to\pn\bm, \quad \pn\rho_n\to\pn\br, \quad\mbox{both weakly in }\,L^2(0,T;H^{1/2}(\Gamma)),
\eeq
and it obviously holds $\,\bm(0)=\mu_0$, $\br(0)=\rho_0$, and $\brg(0)=\rho_{0_{|\Gamma}}$. 
In addition, it is easily verified that 
\beq\label{con8}
\mu_n\,g'(\rho_n)\to \bm\,g'(\br), \quad \mu_n\,g'(\rho_n)\,\pt\rho_n\to \bm\,g'(\br)\,\pt\br, \quad \mbox{both weakly in }\,L^2(Q).
\eeq 
Now, we let $n\to\infty$ in the system \eqref{ss1}--\eqref{ss5}, written for
$\,(\mu_n,\rho_n, \rho_{n_\Gamma})\,$ and the right-hand side $u_{\Gamma,n}$. It then follows
from the above convergence results that $(\bm,\br,\brg)$ solves \eqref{ss1}-\eqref{ss5}
with the \rhs~\,$\bug$, that is, we have $\,(\bm,\br,\brg)=\sol(\bug)$, whence we infer that
the pair $\,((\bm,\br,\brg),\bug)\,$ is admissible for {\bf (CP)}. Its optimality is then
a simple consequence of the weak sequential semicontinuity properties of the cost functional
$\CJ$.\qed

\vspace{4mm}
We now turn our interest to the derivation of first-order necessary optimality conditions
for problem {\bf (CP)}. For this purpose, we generally assume that the hypotheses (A1)--(A6)
are fulfilled and that $\,\bug\in\uad\,$ is an optimal control with associated state
$(\bm,\br,\brg)=\sol(\bug)$ having the properties \eqref{regmu}--\eqref{regrg} and \eqref{separ}.
We aim to eliminate the quantities \,$\eta,\zeta,\zeta_\Gamma\,$ from the variational 
inequality \eqref{vug1}. To this end, we invoke the adjoint state system associated with
\eqref{ss1}--\eqref{ss5} for $\,\bug$, which is formally given by:
\begin{align}
\label{as1}
&-(1+2g(\br))\,p_t-g'(\br)\,\br_t\,p-\Delta p\,=\,g'(\br)\,q+\beta_1(\bm-\hat\mu_Q)\quad
\mbox{in }\,Q,\\[2mm]
\label{as2}
&\pn p=0\quad\mbox{on }\,\Sigma,\quad p(T)=0\quad\mbox{in }\,\oma,\\[2mm]
\label{as3}
&-q_t-\Delta q\,+\,(f''(\br)+\pi'(\br)-\bm\,g''(\br))\,q\non\\
&\quad {}=g'(\br)\left(\bm\,p_t-\bm_t\,p\right)+\beta_2(\br-\hat \rho_Q) \quad
\mbox{in }\,Q,\\[2mm]
\label{as4}
&\pn q-\pt q_\Gamma-\Delta_\Gamma q_\Gamma+(f''_\Gamma(\brg)+\pigs(\brg))\,q_\Gamma
\,=\,\beta_3(\brg-\hat\rho_\Sigma), \quad q_\Gamma=q\onSigma, \quad\mbox{on }\,\Sigma,\\[2mm]
\label{as5}
&q(T)=\beta_4(\br(T)-\hat\rho_\oma)\quad\mbox{in }\,\oma, \quad
q_\Gamma(T)=\beta_5(\brg(T)-\hat\rho_\Gamma)\quad\mbox{on }\,\Gamma.
\end{align}

At this point, we simplify the problem somewhat by imposing the following additional condition:

\vspace{3mm}\noindent
(A7)\quad\,It holds that \revis{$\,\,\beta_4(\br(T)-\hat\rho_\oma) \in {\cal V}\,  \, $ and $\,\, \beta_4(\br(T)-\hat\rho_\oma)_{|\Gamma} = \beta_5(\brg(T)-\hat\rho_\Gamma)\, $ on $\Gamma$}.

\vspace{3mm}\noindent
Observe that (A7) is \revis{satisfied if $\beta_4=\beta_5$, $\hat\rho_\oma
\in V$, $\hat \rho_\Gamma\in\Vg$, and $\,\hat\rho_\Gamma=\hat \rho_{{\oma}_{|\Gamma}}$. In view of the 
fact that always $\,\br(T)\in {\cal V}$, these conditions for the target functions
$\,\hat\rho_\oma\,$ and $\,\hat\rho_\Gamma\,$ seem to be quite natural. A trivial situation is of course obtained when $\beta_4=\beta_5=0$.}

We have the following result.

\vspace{5mm}\noindent
{\sc Theorem 4.2:} \quad\,{\em Suppose that} (A1)--(A6) {\em hold true and that $\bug\in\uad$
is an optimal control whose associated state $(\bm,\br,\brg)=\sol(\bug)$ fulfills} (A7). {\em Then the
adjoint state system} \eqref{as1}--\eqref{as5} {\em has a unique solution $(p,q,q_\Gamma)$ such that}
\begin{align}
\label{regp}
&p\in H^1(0,T;H)\cap C^0([0,T];V)\cap L^2(0,T;W),\\
\label{regq}
&q\in H^1(0,T;H)\cap C^0([0,T];V)\cap L^2(0,T;\Hdue),\\
\label{regqg}
&q_\Gamma\in H^1(0,T;\Hg)\cap C^0([0,T];\Vg)\cap L^2(0,T;H^2(\Gamma)).
\end{align}

\vspace{3mm}\noindent
{\sc Proof:}  \quad\,First, we rewrite the backward-in-time system \eqref{as1}--\eqref{as5}. To this
end, we define the functions
\begin{align*}
&\mtil(x,t):=\bm(x,T-t),\quad\rtil(x,t):=\br(x,T-t),\quad\rgtil(x,t):=\br_\Gamma(x,T-t),
\\[1mm]
&\mqtil(x,t):=\hat\mu_Q(x,T-t),\quad \rqtil:=\hat\rho_Q(x,T-t),\quad
\rstil(x,t):=\hat\rho_\Sigma(x,T-t),
\end{align*}
and consider the initial-boundary value problem
\begin{align}
\label{yztil1}
&(1+2g(\rtil))\,\pt y + g'(\rtil)\,\pt\rtil\,y - \Delta y\,=\,g'(\rtil)\,z
+\beta_1(\mtil-\mqtil)\quad\aeQ,\\[1mm]
\label{yztil2}
&\pn y=0\quad\aeS,\quad y(0)=0\quad\aeO,
\end{align}
\begin{align}
\label{yztil3}
&\pt z-\Delta z+(f''(\rtil)+\pi'(\rtil)-\mtil\,g''(\rtil))\,z\non\\
&\quad {}=g'(\rtil)(\pt\mtil\,y-\mtil\,\pt y)+\beta_2(\rtil-\rqtil)\quad\aeQ,\\[1mm]
\label{yztil4}
&\pn z+\pt \zg-\delg\zg+(f_\Gamma''(\rgtil)+\pigs(\rgtil))\,\zg\,=\,
\beta_3(\rgtil-\rstil)\non\\
&\gianni{\aand\zg=z\onSigma}, \quad\aeS,\\[1mm]
\label{yztil5}
&z(0)=\beta_4(\rtil(0)-\hat\rho_\oma)\quad\aeO,\quad \zg(0)=\beta_5(\rgtil(0)
-\hat\rho_\Gamma)\quad\aeG.
\end{align}
Obviously, any sufficiently smooth solution $(y,z,\zg)$ to \eqref{yztil1}--\eqref{yztil5}
induces a solution $(p,q,q_\Gamma)$ to the adjoint system \eqref{as1}--\eqref{as5}
(and vice versa) by putting
\beq
\label{yzpq}
p(x,t):=y(x,T-t), \quad q(x,t):=z(x,T-t),\quad q_\Gamma(x,t)=\zg(x,T-t).
\eeq
Observe that, thanks to assumption (A7), we have 
\revis{$z(0)\in {\cal V}$}. In addition, we recall the global bounds
\eqref{ssbounds1}, \eqref{ssbounds2}, and the regularity result \eqref{addregmu},
which yield, in particular, that
\begin{align}
\label{yztil6}
&a:=f''(\rtil)+\pi'(\rtil)-\mtil\,g''(\rtil)\in L^\infty(Q),\quad
a_\Gamma:=f_\Gamma''(\rgtil)+\pigs(\rgtil)\in L^\infty(\Sigma),\non\\[1mm]
&\pt\mtil\in L^2(0,T;L^6(\oma)),\quad\pt\rtil\in L^2(0,T;V). 
\end{align}
We aim to show that the system \eqref{yztil1}--\eqref{yztil5} has a unique solution triple
$(y,z,\zg)$ having the same regularity as requested for $(p,q,q_\Gamma)$ in
\eqref{regp}--\eqref{regqg}. We divide the proof of this claim into several steps.

\vspace{5mm}\noindent
\underline{\sc Step 1:}
\par\nobreak
\vspace{1mm}\noindent
We first prove uniqueness. To this end, suppose that two solutions
$(y_i,z_i, z_{i_\Gamma})$, $i=1,2$,
with the asserted regularity are given. Then the triple $(y,z,\zg)$, where $y:=y_1-y_2$,  $z:=z_1-z_2$,
$\zg:=z_{1_\Gamma}-z_{2_\Gamma}$, satisfies the system that results if in \eqref{yztil1}--\eqref{yztil5} the
terms containing the factors $\,\beta_i$, $1\le i\le 5$, are omitted. In particular, $z(0)=0$ and
$\zg(0)=0$. We then can infer from Lemma 2.3 that, for every $t\in (0,T]$,
\begin{align}
\label{uni1}
&\|z\|^2_{H^1(0,t;H)\cap C^0([0,t];V)\cap L^2(0,t;\Hdue)}\,+\,
\|\zg\|^2_{H^1(0,t;\Hg)\cap C^0([0,t];\Vg)\cap L^2(0,t;H^2(\Gamma))}\non\\
&\le C_1\,\|\sigma\|^2_{L^2(Q_t)},\quad\mbox{with $\,\sigma:=g'(\rtil)\,
(\pt\mtil\,y-\mtil\,\pt y)$},
\end{align} 
where, here and in the remainder of the uniqueness proof, we denote by $C_i$, $i\in \pier{\mathbb{N}}$, positive
constants that depend only on the data of the system and on norms of the solutions. Now, by H\"older's
and Young's inequalities, \pier{we have that}
\begin{align}
\label{uni2}
\|\sigma\|^2_{L^2(Q_t)}\,&\le\,C_2\txinto|\pt\mtil\,y-\mtil\,\pt y|^2\dx\ds\non\\
&\le\,C_3\txinto|\pt y|^2\dx\ds\,+\,C_4\tint\|\pt\mtil(s)\|_{\juerg 4}^2\,\|y(s)\|_{\juerg 4}^2\ds\non\\
&\le\,C_3\txinto|\pt y|^2\dx\ds\,+\,C_5\tint\|\pt\mtil(s)\|_6^2\,\|y(s)\|_V^2\ds,
\end{align}
where the mapping $\,s\mapsto \|\pt\mtil(s)\|_6^2\,$ belongs to $L^1(0,T)$.

Next, we add $\,y\,$ on both sides of the equation resulting from \eqref{yztil1},
multiply by $\,\pt y$, and integrate over $Q_t$, where 
$t\in (0,T]$. Since $\,g(\rtil)\ge 0$, we obtain from H\"older's and Young's
inequalities that
\begin{align}
\label{uni3}
&\txinto|\pt y|^2\dx\ds\,+\,\frac 12\,\|y(t)\|_V^2\,\le C_6
\txinto|\pt y|\left(|z|\,+\,(1+|\pt\rtil|)\,|y|\right)\dx\ds\non\\
&\le\,\frac 14\txinto\!|\pt y|^2\dx\ds\,+\,C_7\txinto\!(y^2+z^2)\dx\ds\,+\,
C_8\tint\!\|\pt\rtil(s)\|_4\,\|y(s)\|_4\,\|\pt y(s)\|_2\ds\non\\
&\le\,\frac 12\txinto\!|\pt y|^2\dx\ds\,+\,C_7\txinto\!(y^2+z^2)\dx\ds\,+\,
C_9\tint\!\|\pt\rtil(s)\|_V^2\,\|y(s)\|_V^2\ds\,,
\end{align}
where the mapping $\,s\mapsto \|\pt\rtil(s)\|_V^2\,$ belongs to $L^1(0,T)$.

Now, we multiply the inequality \eqref{uni3} by \gianni{$\,4\,C_1\,C_3$}\, and add the result
to the inequality \eqref{uni1}. Taking \eqref{uni2} into account, we then
conclude from Gronwall's lemma that $\,y=z=0\,$ in $Q$ and $\zg=0$ on
$\Sigma$, \pier{whence} the uniqueness is proved.

\vspace{5mm}\noindent
\underline{Step 2:} 
\par\nobreak
\vspace{1mm}\noindent
We now approximate the system \eqref{yztil1}--\eqref{yztil5}, where we employ
a \juerg{similar} approach as in the proof of Theorem 3.1. 
To this end, we consider for $\tau=\gianni{\tau_N}:=T/N$,
$N\in\nz$, the retarded system
\begin{align}
\label{yztau1}
&(1+2g(\rtil))\,\pt y^\tau + g'(\rtil)\,\pt\rtil\,\TT(y^\tau) - \Delta y^\tau\non\\
&\quad =\,g'(\rtil)\,\TT(z^\tau)+\beta_1(\mtil-\mqtil)\quad\aeQ,\\[1mm]
\label{yztau2}
&\pn y^\tau=0\quad\aeS,\quad y^\tau(0)=0\quad\aeO,\\[1mm]
\label{yztau3}
&\pt z^\tau-\Delta z^\tau+a\,z^\tau\,=\,
g'(\rtil)(\pt\mtil\,y^\tau-\mtil\,\pt y^\tau)+\beta_2(\rtil-\rqtil)\quad\aeQ,\\[1mm]
\label{yztau4}
&\pn z^\tau+\pt z_\Gamma^\tau-\delg z_\Gamma^\tau
+a_\Gamma\,z_\Gamma^\tau\,=\,
\beta_3(\rgtil-\rstil), \quad z^\tau\onSigma=z_\Gamma^\tau \quad\aeS,\\[1mm]
\label{yztau5}
&z^\tau(0)=\beta_4(\rtil(0)-\hat\rho_\oma)\quad\aeO,\quad z_\Gamma^\tau(0)=\beta_5(\rgtil(0)
-\hat\rho_\Gamma)\quad\aeG,
\end{align}
with the translation operator $\TT$ introduced in \eqref{defTT}, and where 
$\,a\,,\,a_\Gamma\,$ are defined in \eqref{yztil6}. Putting again
$\taun:=T/N$, $t_n:=n\taun$, and $I_n:=(0,t_n)$, $1\le n\le N$, for fixed
$N\in\nz$, we then consider for every $n\in\{1,\ldots,N\}$ the initial-boundary value
problem
\begin{align}
\label{yztn1}
& (1+2g(\rtil))\,\pt y_n + g'(\rtil)\,\pt\rtil\,\TTN(y_{n-1}) - \Delta y_n \non \\[1mm]
&\quad \pier{=
\,g'(\rtil)\,\TTN(z_{n-1})
+\beta_1(\mtil-\mqtil) \quad
\mbox{a.\,e. in }\,\oma\times I_n,}\\[1mm]
\label{yztn2}
&\pn y_n=0\quad\mbox{a.\,e. on}\,\Gamma\times I_n,\quad y_n(0)=0\quad\mbox{a.\,e.
in }\,\Omega,  \\[1mm]
\noalign{\allowbreak}
\label{yztn3}
&\pt z_n-\Delta z_n+a\,z_n\,=\,
g'(\rtil)(\pt\mtil\,y_n-\mtil\,\pt y_n)+\beta_2(\rtil-\rqtil)\quad
\mbox{a.\,e. in }\,\oma\times I_n,  \\[1mm]
\label{yztn4}
&\pn z_n+\pt z_{n_\Gamma}-\delg z_{n_\Gamma}
+a_\Gamma\,z_{n_\Gamma}\,=\,
\beta_3(\rgtil-\rstil), \quad z_{n\onSigma}=z_{n_\Gamma} \quad
\mbox{a.\,e. on }\,\Gamma\times I_n,  \\[1mm]
\label{yztn5}
& z_n(0)=\beta_4(\rtil(0)-\hat\rho_\oma)\quad\aeO,\quad z_{n_\Gamma}(0)=\beta_5(\rgtil(0)-\hat\rho_\Gamma)\quad\aeG.
\end{align}
Here, it is understood that $\TTN(y_{n-1})=0$ and $\TTN(z_{n-1})=
\beta_4(\rtil(0)-\hat\rho_\oma)$ for $n=1$.
Using induction with respect to $n$, we again find that \eqref{yztn1}--\eqref{yztn5}
has for every $n\in\{1,\ldots, N\}$ a unique solution with the requested regularity.
Once more, we confine ourselves to show the induction step $n-1\longrightarrow n$. So\pier{,} let
$1<n\le N$, and assume that  for $1\le k\le n-1$ \pier{the} unique solutions 
$(y_k,z_k,z_{k_\Gamma})$ have already been constructed that satisfy the conditions
\begin{align}
\label{glatt}
&y_k\in H^1(I_k;H)\cap C^0(\bar I_k;V)\cap L^2(I_k;W),\non\\[1mm]
&z_k\in H^1(I_k;H)\cap C^0(\bar I_k;V)\cap L^2(I_k;\Hdue),\non\\[1mm]
&z_{k_\Gamma}\in H^1(I_k;\Hg)\cap C^0(\bar I_k;\Vg)\cap L^2(I_k;H^2(\Gamma)).
\end{align}
Since $\,\rtil\in \CQ$ and $\pt\rtil\in L^2(0,T;V)$, we obviously have that 
$$1+2g(\rtil)\in\CQ, \quad \,g'(\rtil)\,\TTN(z_{n-1})
-g'(\rtil)\,\pt\rtil\,\TTN(y_{n-1})+\beta_1(\mtil-\mqtil)\in L^2(I_n;H).
$$
We thus can infer from, e.\,g., \cite[Thm.~2.1]{DHP} that the initial-boundary
value problem \eqref{yztn1}--\eqref{yztn2} enjoys a unique solution
$y_n\in H^1(I_n;H)\cap C^0(\bar I_n;V)\cap L^2(I_n;W)$.  
We then substitute $y_n$ in \eqref{yztn3}, recalling that \revis{$z_n(0)
\in {\cal V}$}. Moreover, we readily verify that
$$
g'(\rtil)(\pt\mtil\,y_n-\mtil\,\pt y_n)+\beta_2(\rtil-\rqtil)\in L^2(I_n;H),\quad
\beta_3(\rgtil-\rstil)\in L^2(I_n;\Hg).
$$ 
Hence, we can infer from Lemma 2.3 the existence of a unique solution pair
$(z_n,z_{n_\Gamma})$ with
\begin{align*}
&z_n\in H^1(I_n;H)\cap C^0(\bar I_n;V)\cap L^2(I_n;\Hdue),\non\\[1mm]
&z_{n_\Gamma}\in H^1(I_n;\Hg)\cap C^0(\bar I_n;\Vg)\cap L^2(I_n;H^2(\Gamma)).
\end{align*}
Arguing as in the proof of Theorem 3.1, we then conclude that $(y_N,z_N,z_{N_\Gamma})$
is the unique solution to the retarded problem \eqref{yztau1}--\eqref{yztau5}
for $\tau=\taun$. 

\vspace{5mm}\noindent
\underline{\sc Step 3:}
\par\nobreak
\vspace{1mm}\noindent
In this part of the existence proof, we derive a priori estimates for the
approximations $(y_N,z_N,z_{N_\Gamma})$, $N\in\nz$, where we denote by 
\pier{$C_i$, $i\in \mathbb{N}$, positive} constants that may depend on the data but not on $N\in\nz$. For the
\revis{sake of better} readability, we omit the superscript $\taun$ in the
estimates, writing it only at the end of each estimation step.    

\vspace{3mm}\noindent
\underline{\sc First estimate:}
\par\nobreak
\vspace{1mm}\noindent
We add $\,y\,$ on both sides of \eqref{yztau1}, multiply the resulting identity
by $\,\pt y$, and integrate over $Q_t$, where $0<t\le T$. Since $g(\rtil)\ge 0$,
we then find that
\beq\label{test1}
\txinto|\pt y|^2\dx\ds\,+\,\frac 12\,\|y(t)\|_V^2\,\le\,\sum_{j=1}^4\,I_j,
\eeq 
where the quantities $I_j$, $1\le j\le 4$, are specified and estimated below.
Clearly, by Young's inequality \pier{and (A6) we infer that}
\begin{align}\label{test2}
&I_1:=\txinto y\,\pt y\dx\ds\,\le\,\frac 15\txinto|\pt y|^2\dx\ds\,+\,C_1
\txinto|y|^2\dx\ds,\\[1mm]
\label{test3}
&I_4:=\beta_1\txinto(\mtil-\mqtil)\,\pt y\dx\ds\,\le\,\frac 15 
\txinto|\pt y|^2\dx\ds\,+\,C_2\,,\\[1mm]
\label{test4}
&\pier{I_3}:=\txinto g'(\rtil)\,\TTN(z)\,\pt y\dx\ds\,\le\,\frac 15
\txinto|\pt y|^2\dx\ds\,+\,C_3\txinto|\TTN(z)|^2\dx\ds\non\\
&\hspace*{9mm}\le\,\frac 15\txinto|\pt y|^2\dx\ds\,+\,C_4\txinto|z|^2\dx\ds\,+\,C_5,
\end{align}
where in the last estimate we have \pier{employed \eqref{TT1}.} Finally, we
argue as in the estimates \eqref{estn9}--\eqref{estn11} to conclude that
\begin{align}
\label{test5}
\pier{I_2}:&=\,-\txinto g'(\rtil)\pt\rtil\,\TTN(y)\,\pt y\dx\ds\non\\
&\le\,
\frac 15\txinto|\pt y|^2\dx\ds\,+\,C_6\tint\psi(s)\|y(s)\|_V^2\ds,
\end{align}
where the function
$$
s\mapsto\psi(s):=\left\{
\begin{array}{ll}
\|\pt\rtil(s+\tau_N)\|_V^2&\mbox{if }\,0\le s\le T-\tau_N\\[1mm]
0&\mbox{if }\,T-\tau_N< s\le T
\end{array}
\right.
$$
is bounded in $L^1(0,T)$, uniformly in $N\in\nz$. Combining \eqref{test1}--\eqref{test5}, we have thus
shown the estimate
\beq\label{test6}
\frac 15\,\|\pt y\|^2_{L^2(Q_t)}\,+\,\frac 12\,\|y(t)\|_V^2
\,\le\, C_7\,+\,C_8\,\|z\|^2_{L^2(Q_t)}\,+\,C_9\tint (1+\psi(s))\|y(s)\|_V^2\ds\,.
\eeq

\vspace{2mm}\noindent
\underline{\sc Second estimate:}
\par\nobreak
\vspace{1mm}\noindent
Next, we add $\,z\,$ on both sides of \eqref{yztau3}, and $\,\zg\,$ on both 
sides of \eqref{yztau4}, and multiply the first resulting equation by $\pt z$.
Integrating over $Q_t$, where $0<t\le T$, we find the inequality
\begin{align}
\label{test7}
&\txinto|\pt z|^2\dx\ds\,+\tginto|\pt\zg|^2\dg\ds\,+\,\frac 12\left(
\|z(t)\|_V^2+\|\zg(t)\|_{\Vg}^2\right)\non\\
&\le\,\txinto\left(1+\|a\|_{L^\infty(Q)}\right)|z|\,|\pt z|\dx\ds\,+
\txinto\left(1+\|a_\Gamma\|_{L^\infty(\Sigma)}\right)|\zg|\,|\pt\zg|\dg\ds\non\\
&\quad \,+\,\beta_2\txinto(\rtil-\rqtil)\,\pt z\dx\ds
\,+\,\beta_3\txinto(\rgtil-\rstil)\,\pt \zg\dg\ds\non\\
&\quad\,+\,C_{10}\txinto\left(|\pt\mtil|\,|y|+|\mtil|\,|\pt y|\right)\,|\pt z|\dx\ds\non\\
&\quad \gianni{{}+\frac 12 \bigl(\|\beta_4(\rtil(0)-\hat\rho_\oma)\|_V^2 + \|\beta_5(\rgtil(0)-\hat\rho_\Gamma)\|_{V_\Gamma}^2 \bigr)}
\end{align}
\gianni{and observe that the terms in the last line are finite by assumption~(A7).}
Thanks to (A6) and Young's inequality, the first four summands on the \rhs~are bounded
by an expression of the form
\beq
\frac 14\left(\|\pt z\|^2_{L^2(Q_t)}+\|\pt\zg\|^2_{L^2(\Sigma_t)}\right)
+\,C_{11}\left(1\,+\,\|z\|^2_{L^2(Q_t)}\,+\,\|\zg\|^2_{L^2(\Sigma_t)}\right)\,.
\eeq
Moreover, since $\mtil\in L^\infty(Q)$, we have
\beq
C_{10}\txinto|\mtil|\,|\pt y|\,|\pt z|\dx\ds\,\le\,\frac 14\|\pt z\|^2_{L^2(Q_t)}\,+\,
C_{12}\,\|\pt y\|^2_{L^2(Q_t)}\,.
\eeq
In addition, by also using H\"older's inequality,
\begin{align}
\label{test8}
&C_{10}\txinto|\pt\mtil|\,|y|\,|\pt z|\dx\ds\,\le\,C_{13}\tint
\|\pt\mtil(s)\|_{\pier{6}}\,\|y(s)\|_{\pier{3}}\,\|\pt z(s)\|_2\ds\non\\
&\le\,\frac 14\|\pt z\|^2_{L^2(Q_t)}\,+\,C_{14}
\tint \|\pt\mtil(s)\|^2_6\,\|y(s)\|_V^2\ds\,.
\end{align}
Combining the estimates \eqref{test7}--\eqref{test8}, we have thus shown that
\begin{align}
\label{test9}
&\frac 14\,\|\pt z\|^2_{L^2(Q_t)}\,+\,\frac 34\,\|\pt\zg\|^2_{L^2(\Sigma_t)}
\,+\,\frac 12\left(
\|z(t)\|_V^2+\|\zg(t)\|_{\Vg}^2\right)\non\\
&\le\,C_{11}\left(1\,+\,\|z\|^2_{L^2(Q_t)}\,+\,\|\zg\|^2_{L^2(\Sigma_t)}\right)
\,+\, C_{12}\,\|\pt y\|_{L^2(Q_t)}^2\non\\
&\quad\,+\,C_{14}
\tint \|\pt\mtil(s)\|^2_6\,\|y(s)\|_V^2\ds
\gianni{{}+C_{15}}\,,
\end{align}
where the mapping $\,s\mapsto \|\pt\mtil(s)\|^2_6\,$ belongs to $L^1(0,T)$.

Now, we multiply \eqref{test6} by $\,10\,C_{12}\,$ and add the resulting
inequality to \eqref{test9}. It then follows from Gronwall's lemma that, for
all $t\in (0,T]$ and $N\in\nz$,
\begin{align}
\label{test10}
&\|y^{\taun}\|_{H^1(0,t;H)\cap L^\infty(0,t;V)}\,+\,
\|z^{\taun}\|_{H^1(0,t;H)\cap L^\infty(0,t;V)}\non\\
&\quad {}+\|\zg^{\taun}\|_{H^1(0,t;\Hg)
\cap L^\infty(0,t;\Vg)}\,\le\, \gianni{C_{16}}\,.
\end{align}

\vspace{2mm}\noindent
\underline{\sc Third estimate:}
\par\nobreak
\vspace{1mm}\noindent
Now that the basic estimate \eqref{test10} is shown, we can easily conclude from
comparison in \eqref{yztau1} and \eqref{yztau3}, respectively, that
\beq\label{test11}
\|\Delta y\|_{L^2(Q)}\,+\,\|\Delta z\|_{L^2(Q)}\,\le\,\gianni{C_{17}},
\eeq
whence, using \pier{the boundary condition in \eqref{yztau2} and} standard elliptic estimates, \pier{we deduce that}
\beq\label{test12}
\|y\|_{L^2(0,T;W)}\,\le\,\gianni{C_{18}}\,.
\eeq
\pier{Moreover}, we invoke \cite[Thm.~3.2, p.~1.79]{brez} to conclude that
$$\int_0^T\|z(t)\|^2_{H^{3/2}(\oma)}\dt\,\le\,\gianni{C_{19}}
\int_0^T\left(\|\Delta z(t)\|^2_H+\|\zg(t)\|_{H^1(\Gamma)}^2\right)\dt,
$$
which entails that
\beq
\|z\|_{L^2(0,T;H^{3/2}(\oma))}\,\le\,\gianni{C_{20}}\,. \label{test12bis}
\eeq
Hence, by the trace theorem (cf. \cite[Thm.~2.27, p.~1.64]{brez}), we have that
\beq\label{test13}
\|\pn z\|_{L^2(0,T;\Hg)}\,\le\,\gianni{C_{21}}\,.
\eeq
Comparison in \eqref{yztau4} then yields that
\beq\label{test14}
\|\delg\zg\|_{L^2(0,T;\Hg)}\,\le\,\gianni{C_{22}}\,,
\eeq
and it follows from the boundary version of the elliptic estimates that
\beq\label{test15}
\|\zg\|_{L^2(0,T;H^2(\Gamma))}\,\le\,\gianni{C_{23}}\,.
\eeq
\pier{At this point, \eqref{test10}, \eqref{test11}, \eqref{test15}
allow us to improve \eqref{test12bis} as}
\beq\label{test16}
\|z\|_{L^2(0,T;\Hdue)}\,\le\,\gianni{C_{24}}\,.
\eeq
Recalling that the embeddings $\,(H^1(0,T;H)\cap L^2(0,T;\Hdue))\subset
C^0([0,T];V)$ and \linebreak $\,(H^1(0,T;\Hg)\cap L^2(0,T;H^2(\Gamma))\subset
C^0([0,T];\Vg)$ are continuous, we have finally shown the estimate
\begin{align}
\label{test17}
&\|y^{\taun}\|_{H^1(0,T;H)\cap C^0([0,T];V)\cap L^2(0,T;W)}\,+\,
\|z^{\taun}\|_{H^1(0,T;H)\cap C^0([0,T];V)\cap L^2(0,T;\Hdue)}\non\\[1mm]
&+\,\|z_\Gamma^{\taun}\|_{H^1(0,T;\Hg)\cap C^0([0,T];\Vg)\cap L^2(0,T;
H^2(\Gamma))}\,\le\,\gianni{C_{25}}\,.
\end{align}

\vspace{3mm}\noindent
\underline{\sc Step 4:}
\par\nobreak
\vspace{1mm}\noindent
We now conclude the existence part of the proof. To this end, we observe that
\eqref{test17} yields the existence of a triple $(y,z,\zg)$ such that, at
least for a subsequence which is again indexed by $N$, \pier{we have 
that}
\begin{align}
\label{conv1}
y^{\tau_N}&\to y\quad\,\,\mbox{weakly in }\,H^1(0,T;H)\cap C^0([0,T];V)
\cap L^2(0,T;W)\,,\\[1mm]
\label{conv2}
z^{\tau_N}&\to z\quad\,\,\mbox{weakly in }\,H^1(0,T;H)\cap C^0([0,T];V)
\cap L^2(0,T;\Hdue)\,,\\[1mm]
\label{conv3}
z_\Gamma^{\tau_N}&\to \zg\quad\mbox{weakly in }\,H^1(0,T;\Hg)\cap C^0([0,T];\Vg)
\cap L^2(0,T;H^2(\Gamma))
\end{align} 
\pier{as $N\to\infty$.} We are now \pier{in} a similar situation as in the proof of Theorem 3.1 after showing the
corresponding convergence results \eqref{coneta}--\eqref{conzetg}. Adapting
the arguments used there (with obvious modifications) to our situation, we can
conclude that $\,(y,z,\zg)\,$ is in fact a solution to the transformed system
\eqref{yztil1}--\eqref{yztil5} having the asserted regularity properties. As this is
a rather straightforward repetition of the argumentation utilized there, we may allow ourselves to leave it to the reader to work out
the details. Since, as it was shown in Step 1, such a solution is uniquely 
determined, we can conclude that the adjoint state system \eqref{as1}--\eqref{as5}
has indeed a unique solution satisfying \eqref{regp}--\eqref{regqg}. 
The assertion is thus completely proved.
\qed

\vspace{5mm}
We now can eliminate the functions $(\eta,\zeta,\zetg)$ from the 
variational inequality \eqref{vug1}. We have the following result.

\vspace{3mm}\noindent
{\sc Corollary 4.3:} \quad\,{\em Suppose that} (A1)--(A6) {\em are satisfied, assume that 
$\bug\in \uad$
is an optimal control whose associated state $(\bm,\br,\brg)=\sol(\bug)$ fulfills} (A7), 
{\em and let $(p,q,q_\Gamma)$ be the corresponding unique solution to the adjoint state system} 
\eqref{as1}--\eqref{as5}
{\em established in Theorem 4.2. Then \pier{there} holds  the variational inequality}
\beq
\label{vug2}
\teginto\left(q_\Gamma+\beta_6\,\bug\right)\left(v_\Gamma-\bug\right)\dg\dt\,\ge\,0\,\quad\forall\,
v_\Gamma\in\uad.
\eeq  

\vspace{3mm}\noindent
{\sc Proof:} \quad\,Let $\,v_\Gamma\in\uad\,$ be arbitrary, and let $\,h:=v_\Gamma-\bug$.
We multiply \eqref{ls1} by~$\,p$, \gianni{\eqref{ls3}} by~$\,q$, \eqref{as1} by~$\,-\eta$, 
\gianni{\eqref{as3}} by~$\,-\zeta$, add the resulting equations and integrate over~$Q$
and by parts. A straightforward calculation, in which many terms cancel out, leads to the
identity
\begin{align}
\label{vug3}
&\texinto\pt\bigl((1+2g(\br))\,p\,\eta+q\,\zeta\bigr)\dx\dt \,+\,\teginto(\zetg\,\pn q
-q_\Gamma\,\pn\zeta)\dg\dt\non\\
&+\texinto \left\{\pt\bm\,g'(\br)\,\zeta\,p+\bm\,g''(\br)\,\pt\br\,\zeta\,p
+\bm\,g'(\br)\,\pt\zeta\,p+\bm\,g'(\br)\,\zeta\,\pt p\right\}\dx\dt\non\\
&=\,-\,\beta_1\texinto(\bm-\hat\mu_Q)\,\eta\dx\dt\,-\,\beta_2\texinto(\br-\hat\rho_Q)\,\zeta\dx\dt\,.
\end{align}
\gianni{Clearly}, the integrand in the curly bracket\pier{s} equals $\,\pt(\bm\,g'(\br)\,\zeta\,p)$,
\gianni{whence the corresponding integral vanishes since $\zeta(0)=p(T)=0$}. 
Moreover, owing to \eqref{ls4} and~\eqref{as4},
\beq
\label{vug4}
\teginto(\zetg\,\pn q-q_\Gamma\,\pn\zeta)\dg\dt\,=\,\teginto\bigl(\pt(q_\Gamma\,\zetg)+\beta_3\,(\brg-\hat\rho_\Sigma)\zetg
-q_\Gamma(v_\Gamma-\bug)\bigr)\dg\dt\,,
\eeq 
because the terms involving the Laplace-Beltrami operator cancel each other.
Recalling that $\eta(0)=\zeta(0)=p(T)=0$ in $\oma$ and $\zetg(0)=0$ on $\Gamma$, invoking
the end point conditions \eqref{as5}, and rearranging terms, we finally arrive at the identity 
\begin{align}
\label{vug5}
&\beta_1\texinto (\bm-\hat\mu_Q)\,\eta\dx\dt\,+\,\beta_2\texinto(\br-\hat\rho_Q)\,\zeta\dx\dt
\,+\,\beta_3\teginto(\br_\Gamma-\hat\rho_\Sigma)\,\zeta_\Gamma\dg\dt
\nonumber\\[1mm]
&+\,\beta_4\xinto(\br(T)-\hat\rho_\oma)\,\zeta(T)\dx\,+\,\beta_5\ginto(\br_\Gamma(T)-\hat\rho_\Gamma)
\,\zeta_\Gamma(T)\dg\nonumber\\[1mm]
&=\,\teginto q_\Gamma\,(v_\Gamma-\bug)\dg\dt\,.
\end{align}
The assertion then follows from insertion of this identity in \eqref{vug1}.
\qed

\vspace{5mm}\noindent
{\sc Remark 4.4:} \quad\,If $\,\beta_6>0$, then \eqref{vug2} implies that $\bug$ is nothing but the $L^2(\Sigma)$-orthogonal projection of $\,-\beta_6^{-1} q_\Gamma\,$ onto $\uad$. \revis{However, recalling (A4) let us notice that $\bug$ is different from the pointwise projection on $[u_*, u^*]$, unless the latter belongs to $\uad$, that is,  it satisfies the additional constraint on the time derivative.} 

{\small
 
}
\end{document}